# HIGH-RESOLUTION ASYMPTOTICS FOR THE ANGULAR BISPECTRUM OF SPHERICAL RANDOM FIELDS[1]

By Domenico Marinucci

*Università di Roma "Tor Vergata"*

In this paper we study the asymptotic behavior of the angular bispectrum of spherical random fields. Here, the asymptotic theory is developed in the framework of fixed-radius fields, which are observed with increasing resolution as the sample size grows. The results we present are then exploited in a set of procedures aimed at testing non-Gaussianity; for these statistics, we are able to show convergence to functionals of standard Brownian motion under the null hypothesis. Analytic results are also presented on the behavior of the tests in the presence of a broad class of non-Gaussian alternatives. The issue of testing for non-Gaussianity on spherical random fields has recently gained enormous empirical importance, especially in connection with the statistical analysis of cosmic microwave background radiation.

**1. Introduction.** Several statistical challenges are now arising in connection with cosmological data, and more precisely, for the analysis of cosmic microwave background radiation (CMB). CMB can be viewed as a snapshot of the Universe approximately $3 \times 10^5$ years after the Big Bang; technological progress has made possible a number of experiments aimed at measuring the properties of this radiation. Pioneering results were released in 1992 by the NASA mission COBE [31], which was the first to release a full-sky map of CMB fluctuations; the statistical properties of these fluctuations were then further investigated by several balloon experiments, starting with BOOMERanG [8] and MAXIMA [14]. A major breakthrough is associated with two satellite missions, namely WMAP [4], which released its first data set in February 2003, with much more detailed data to come in the next four years, and Planck, which is scheduled to be launched in Spring 2007 and expected to provide maps with much greater resolution. Over the next

Received March 2004; revised March 2005.
[1]Supported by MIUR.

ten years, an immense amount of cosmological information is expected from these huge data sets; statistical efforts needed to extract this information are equally challenging and impressive.

From the mathematical point of view, CMB can be represented as a random field $T(\theta,\varphi)$ indexed by the unit sphere $S^2$, that is, for each azimuth $0 \le \theta \le \pi$ and elongation $0 \le \varphi < 2\pi$, $T(\theta,\varphi)$ is a real random variable defined on some probability space. We shall always assume that $T(\theta,\varphi)$ is a zero-mean, finite-variance, mean-square continuous and isotropic random field, that is, its distribution is invariant with respect to the group of rotations. Until very recently, the assumption of an isotropic Universe has been taken for granted in cosmological physics, as a consequence of Einstein's *Cosmological Principle* that the Universe should "appear the same" to an observer located anywhere in space. Quite interestingly, the first release of WMAP has raised some doubts on this condition [10, 15, 27]. Testing for isotropy is a very interesting topic, almost completely open for statistical research; we consider this issue, however, beyond the purpose of the present work.

If an isotropic field is Gaussian, its dependence structure is completely identified by the angular correlation function and its harmonic transform, the angular power spectrum (to be defined in the next section). For non-Gaussian fields, the dependence structure becomes much richer, and higher-order correlation functions are of interest; in turn, this leads to the analysis of so-called higher-order angular power spectra. Because these angular power spectra are identically zero for Gaussian fields, they also provide natural tools to test for non-Gaussianity: this is a topic of greatest importance in modern cosmological data analysis. Indeed, on one hand the validation of the Gaussian assumption is urged by the necessity to provide a firm basis for statistical inference on cosmological parameters, which is dominated by likelihood approaches. More importantly, tests for Gaussianity are needed to discriminate among competing scenarios for the physics of the primordial epochs: here, the currently favored inflationary models predict (very close to) Gaussian CMB fluctuations, whereas other models yield different observational consequences [3, 28, 29]. Tests for non-Gaussianity are also powerful tools to detect systematic effects in the outcome of the experiments. For these reasons, many papers have focused on testing for non-Gaussianity on CMB, some of them by means of topological properties of Gaussian fields (e.g., [9, 13, 26, 36, 37]), others through spherical wavelets [2], and still others by harmonic space methods (e.g., [16, 20, 21, 24]). A short survey of the literature on testing for non-Gaussianity on CMB is in [23].

In this paper, we investigate the asymptotic properties for the observed bispectrum of spherical Gaussian fields, and we analyze its use as a probe of non-Gaussian features. The bispectrum (defined in Section 2) is probably the single most popular statistic to search for non-Gaussianity in CMB data; on one hand, in fact, working on harmonic space is extremely convenient, and the bispectrum is the simplest harmonic space statistic which is sensitive to non-Gaussian features. On the other hand, it is possible to derive analytically the behavior of the bispectrum for non-Gaussian fields of physical interest. In fact, although the procedures considered in the present work are new (to the best of our knowledge), a number of insightful and important papers have already considered the bispectrum for CMB data analysis, for instance, [19, 20, 21, 30]. Few analytical results, however, have so far been produced on the statistical properties of these procedures. The reason for this can be partially explained as follows: as described in the next section, the bispectrum is a function of some spherical harmonic coefficients $a_{lm}$'s; in the presence of an ideal experiment, the latter are easily derived from a map of CMB fluctuations by a harmonic transform performed on the observed data. However, it is important to stress that in realistic situations the $a_{lm}$'s are observed with error, due to instrumental noise, gaps in the maps and many other sources (these problems have been described in [33]). Determining the properties of the procedures in realistic situations that take into account all the features of real life experiments is extremely difficult; most work in CMB is thus based on comparisons of estimates based on real data with the expected results of simulations under a particular null hypothesis. In this paper, we assume that the $a_{lm}$'s are observed without error; this is a simplifying assumption, which we adopt because it seems important to narrow the gap between data analysis practice and its mathematical foundations, at least in an idealized case. Future work, however, should be directed at relaxing this assumption.

In Section 2 we define the bispectrum, taking particular care to discuss the conditions to ensure that it represents a rotationally invariant statistic. The asymptotic behavior of its higher-order moments is considered in Section 3: these results are derived under Gaussianity, but the general technique by which they are established (which adopts a formalism from graph theory) may have some independent interest under broader assumptions. Section 4 considers the effect of an unknown angular power spectrum, whereas Section 5 discusses statistical applications, with asymptotic results in Gaussian and non-Gaussian circumstances. The results we present in this section suggest that consistent tests of Gaussianity can exist even for random fields defined on a bounded domain, which is to some extent unexpected. Section 6 discusses the rationale behind the approach presented and draws some conclusions; some technical results are collected in the Appendix.

**2. The angular bispectrum.** The Fourier transform on the sphere is defined by the spherical harmonics, which can be written explicitly as

$$Y_{lm}(\theta,\varphi)=\sqrt{\frac{2l+1}{4\pi}\frac{(l-m)!}{(l+m)!}}P_{lm}(\cos\theta)\exp(im\varphi)\qquad \text{for } m\geq 0,$$

$$Y_{lm}(\theta,\varphi) = (-1)^m Y^*_{l,-m}(\theta,\varphi) \qquad \text{for } m<0,$$

where the asterisk denotes complex conjugation and $P_{lm}(\cos\theta)$ denotes the associated Legendre polynomial of degree $l,m$, that is,

$$P_{lm}(x) = (-1)^m(1-x^2)^{m/2}\frac{d^m}{dx^m}P_l(x), \qquad P_l(x) = \frac{1}{2^l l!}\frac{d^l}{dx^l}(x^2-1)^l,$$
$$m = 0,1,\dots,l, l = 1,2,3,\dots.$$

A detailed discussion of the properties of the spherical harmonics can be found in ([34], Chapter 5), or in [35]. For isotropic fields, the following spectral representation holds in the mean-square sense (see also [1, 22, 38]):

$$T(\theta,\varphi) = \sum_{l=1}^{\infty}\sum_{m=-l}^{l} a_{lm}Y_{lm}(\theta,\varphi), \tag{1}$$

where the triangular array $\{a_{lm}\}$ represents a set of random coefficients, which can be obtained from $T(\theta,\varphi)$ through the inversion formula

$$a_{lm} = \int_{-\pi}^{\pi}\int_0^{\pi} T(\theta,\varphi)Y^*_{lm}(\theta,\varphi)\sin\theta\, d\theta\, d\varphi, \tag{2}$$
$$m = 0,\pm1,\dots,\pm l, l = 1,2,\dots.$$

These coefficients are complex-valued, zero-mean and uncorrelated; hence, if $T(\theta,\varphi)$ is Gaussian, they have a complex Gaussian distribution, and they are independent over $l$ and $m \geq 0$ [although $a_{l,-m} = (-1)^m a^*_{lm}$], with variance $E|a_{lm}|^2 = C_l$, $m = 0,\pm1,\dots,\pm l$. The index $l$ is usually labeled a multipole and in principle it runs from 1 to infinity; each multipole corresponds approximately to an angular resolution of $180/l°$. In any realistic experiment, however, there is an upper limit (which we denote by $L$) on the multipoles we may observe, depending upon the resolution of the experiment and the presence of noise; $L$ is reckoned to be of the order of 600/800 for WMAP and 2000/2500 for Planck. Strictly speaking, in CMB cosmology the spectral representation (1) is really only defined for $l \geq 2$; the so-called dipole $l = 1$ is in fact dominated by kinematic effects and thus it is removed from the data.

The sequence $\{C_l\}$ denotes the angular power spectrum: we shall always assume that $C_l$ is strictly positive, for all values of $l$. This condition is very mild; to draw an analogy with the theory of stationary random fields defined on $R^d$, it is equivalent to the (very common) assumption that their spectral density is strictly positive at all frequencies. As discussed in Section 5, the condition is met by virtually all models of cosmological interest. A natural estimator for $C_l$ is

$$\widehat{C}_l = \frac{1}{2l+1}\sum_{m=-l}^{l}|a_{lm}|^2, \qquad l = 1,2,\dots, \tag{3}$$

which is clearly unbiased; see also [12]. As mentioned in the Introduction, if the field is Gaussian, the angular power spectrum completely identifies its dependence structure. For non-Gaussian fields, the dependence structure becomes much richer, and higher-order moments of the $a_{lm}$'s are of interest; this leads to the analysis of so-called higher-order angular power spectra.

Generally speaking, the angular bispectrum can be viewed as the harmonic transform of the three-point angular correlation function, much as the angular power spectrum is the Legendre transform of the (two-point) angular correlation function. More precisely, write $\Omega_i = (\theta_i, \varphi_i)$, for $i = 1, 2, 3$; we have

$$\begin{aligned} &ET(\Omega_1)T(\Omega_2)T(\Omega_3) \\ &\quad = \sum_{l_1,l_2,l_3=1}^{\infty} \sum_{m_1,m_2,m_3} B_{l_1 l_2 l_3}^{m_1 m_2 m_3} Y_{l_1 m_1}(\Omega_1) Y_{l_2 m_2}(\Omega_2) Y_{l_3 m_3}(\Omega_3), \end{aligned} \tag{4}$$

where the bispectrum $B_{l_1 l_2 l_3}^{m_1 m_2 m_3}$ is given by

$$B_{l_1 l_2 l_3}^{m_1 m_2 m_3} = E(a_{l_1 m_1} a_{l_2 m_2} a_{l_3 m_3}). \tag{5}$$

Here, and in the sequel, the sums over $m_i$ run from $-l_i$ to $l_i$, unless otherwise indicated. Both (4) and (5) are clearly equal to zero for zero-mean Gaussian fields. Moreover, the assumption that the CMB random field is statistically isotropic entails that the right- and left-hand sides of (4) should be left unaltered by a rotation of the coordinate system. Therefore $B_{l_1 l_2 l_3}^{m_1 m_2 m_3}$ must take values ensuring that the three-point correlation function on the left-hand side of (4) remains unchanged if the three directions $\Omega_1, \Omega_2$ and $\Omega_3$ are rotated by the same angle. Careful choices of the orientations entail that the angular bispectrum of an isotropic field can be nonzero only if:

(a) $l_1, l_2$ and $l_3$ satisfy the triangle rule, $l_i \le l_j + l_k$ for all choices of $i, j, k = 1, 2, 3$,
(b) $l_1 + l_2 + l_3 =$ even and
(c) $m_1 + m_2 + m_3 = 0$.

More generally, Hu [17] shows that a necessary and sufficient condition for $B_{l_1 l_2 l_3}^{m_1 m_2 m_3}$ to represent the angular bispectrum of an isotropic random field is that there exist a real symmetric function of $l_1, l_2, l_3$, which we denote $b_{l_1 l_2 l_3}$, such that we have the identity

$$B_{l_1 l_2 l_3}^{m_1 m_2 m_3} = \mathcal{G}_{l_1 l_2 l_3}^{m_1 m_2 m_3} b_{l_1 l_2 l_3}; \tag{6}$$

$b_{l_1 l_2 l_3}$ is labeled the reduced bispectrum. In (6) we are using the Gaunt integral $\mathcal{G}_{l_1 l_2 l_3}^{m_1 m_2 m_3}$, defined by

$$\begin{aligned} \mathcal{G}_{l_1 l_2 l_3}^{m_1 m_2 m_3} &= \int_0^{\pi} \int_0^{2\pi} Y_{l_1 m_1}(\theta, \varphi) Y_{l_2 m_2}(\theta, \varphi) Y_{l_3 m_3}(\theta, \varphi) \sin\theta \, d\varphi \, d\theta \\ &= \left( \frac{(2l_1+1)(2l_2+1)(2l_3+1)}{4\pi} \right)^{1/2} \begin{pmatrix} l_1 & l_2 & l_3 \\ 0 & 0 & 0 \end{pmatrix} \begin{pmatrix} l_1 & l_2 & l_3 \\ m_1 & m_2 & m_3 \end{pmatrix}, \end{aligned}$$

where the so-called "Wigner $3j$ symbols" appearing on the second line are defined in the Appendix. It can be shown that the Gaunt integral is identically equal to zero unless the conditions (a)–(c) are fulfilled. Often the dependence on $m_1, m_2, m_3$, which does not carry any physical information if the field is isotropic, is eliminated by focusing on the angular averaged bispectrum, defined by

$$\begin{aligned} B_{l_1 l_2 l_3} &= \sum_{m_1=-l_1}^{l_1} \sum_{m_2=-l_2}^{l_2} \sum_{m_3=-l_3}^{l_3} \begin{pmatrix} l_1 & l_2 & l_3 \\ m_1 & m_2 & m_3 \end{pmatrix} B_{l_1 l_2 l_3}^{m_1 m_2 m_3} \\ &= \left( \frac{(2l_1+1)(2l_2+1)(2l_3+1)}{4\pi} \right)^{1/2} \begin{pmatrix} l_1 & l_2 & l_3 \\ 0 & 0 & 0 \end{pmatrix} b_{l_1 l_2 l_3}, \end{aligned} \tag{7}$$

where we have used (48) (see the Appendix). In practice, of course, the bispectrum is not observable; its minimum mean-square error estimator is provided by Hu [17],

$$\widehat{B}_{l_1 l_2 l_3} = \sum_{m_1=-l_1}^{l_1} \sum_{m_2=-l_2}^{l_2} \sum_{m_3=-l_3}^{l_3} \begin{pmatrix} l_1 & l_2 & l_3 \\ m_1 & m_2 & m_3 \end{pmatrix} (a_{l_1 m_1} a_{l_2 m_2} a_{l_3 m_3}).$$

The statistic $\widehat{B}_{l_1 l_2 l_3}$ is called the (sample) angle averaged bispectrum; for any realization of the random field $T$, it is a real-valued scalar, which does not depend on the choice of the coordinate axes, and it is invariant with respect to permutation of its arguments $l_1, l_2, l_3$.

Now note that, under Gaussianity, the distribution of $a_{lm}/C_l^{1/2}$ does not depend on any nuisance parameter. The bispectrum can hence be easily made model-independent; namely, we can focus on the normalized bispectrum, which we define by

$$I_{l_1 l_2 l_3} = (-1)^{(l_1+l_2+l_3)/2} \frac{\widehat{B}_{l_1 l_2 l_3}}{\sqrt{C_{l_1} C_{l_2} C_{l_3}}}. \tag{8}$$

The factor $(-1)^{(l_1+l_2+l_3)/2}$ is usually not included in the definition of the normalized bispectrum; it corresponds, however, to the sign of Wigner's coefficients for $m_1 = m_2 = m_3 = 0$, and thus it seems natural to include it to ensure that $I_{l_1 l_2 l_3}$ and $b_{l_1 l_2 l_3}$ share the same parity [see (7)].

In practice, of course, $I_{l_1 l_2 l_3}$ is infeasible, and it is thus replaced by the statistic

$$\widehat{I}_{l_1 l_2 l_3} = (-1)^{(l_1+l_2+l_3)/2} \frac{\widehat{B}_{l_1 l_2 l_3}}{\sqrt{\widehat{C}_{l_1} \widehat{C}_{l_2} \widehat{C}_{l_3}}};$$

see [12].

**3. Higher-order moments of the angular bispectrum.** In this section we shall investigate the behavior of the higher-order moments for the normalized bispectrum (8), under the assumption of Gaussianity. We assume without loss of generality $l_1 \leq l_2 \leq l_3$, and define

$$\Delta_{l_1 l_2 l_3} \stackrel{\text{def}}{=} 1 + \delta_{l_1}^{l_2} + \delta_{l_2}^{l_3} + 3\delta_{l_1}^{l_3} = \begin{cases} 1, & \text{for } l_1 < l_2 < l_3, \\ 2, & \text{for } l_1 = l_2 < l_3 \text{ or } l_1 < l_2 = l_3, \\ 6, & \text{for } l_1 = l_2 = l_3; \end{cases}$$

here and in the sequel, $\delta_a^b$ denotes Kronecker's delta, that is, $\delta_a^b = 1$ for $a = b$, zero otherwise.

Under Gaussianity, it is obvious that the expectation of all odd powers of $I_{l_1 l_2 l_3}$ is zero. To analyze the behavior of even powers, we first recall that, for a multivariate zero-mean Gaussian vector $(x_1, \ldots, x_{2k})$, we have the following diagram formula:

$$E(x_1 \times x_2 \times \cdots \times x_{2k}) = \sum (E x_{i_1} x_{i_2}) \times \cdots \times (E x_{i_{2k-1}} x_{i_{2k}}), \tag{9}$$

where the sum is over all the $(2k)!/(k!2^k)$ different ways of grouping $(x_1, \ldots, x_{2k})$ into pairs (e.g., [1], page 108). Even powers of $I_{l_1 l_2 l_3}$ yield even powers of the $a_{lm}$'s, which have a complex Gaussian distribution, weighted by Wigner's $3j$ coefficients; we shall then need to use some arguments from graph theory, which is widely used in physics when handling Wigner's $3j$ coefficients (see [34], Chapter 11).

Consider the Cartesian product $I \otimes J$, where $I, J$ are sets of positive integers of cardinality $\#(I) = P, \#(J) = Q$; it is convenient to visualize these elements in a $P \times Q$ matrix with $P$ rows labeled by $i$ and $Q$ columns labeled by $j$. A *diagram* $\gamma$ is any partition of the $P \times Q$ elements into pairs like $\{(i_1, j_1), (i_2, j_2)\}$; these pairs are called the *edges* of the diagram. For our purposes, it is enough to consider diagrams with an even number of rows $P$; we label $\Gamma(I, J)$ the family of these diagrams. It can be checked that, for given $I, J$, there exist $(P \times Q - 1)!!$ different diagrams, each of them composed of $(P \times Q)/2$ pairs; we recall that $(2p-1)!! \stackrel{\text{def}}{=} (2p-1) \times (2p-3) \times \cdots \times 1$ for $p = 1, 2, \ldots$. We also note that if we identify each row $i_k$ with a *vertex* (or *node*), and view these vertices as linked together by the edges $\{(i_k, j_k), (i_{k'}, j_{k'})\} = i_k i_{k'}$, then it is possible to associate to each diagram a *graph*. As it is well known, a graph is an ordered pair $(I, E)$, where $I$ is non-empty (in our case the set of the rows of the diagram), and $E$ is a set of unordered pairs of vertices (in our case, the pairs of rows that are linked in a diagram). We consider only graphs which are not directed, that is, $(i_1 i_2)$ and $(i_2 i_1)$ identify the same edge; however, we do allow for repetitions of edges (two rows may be linked twice), in which case the term *multigraph* is more appropriate. In general, a graph carries less information than a diagram (the information on the "columns," i.e., the second element $j_k$, is neglected), but it is much easier to represent pictorially. We shall use some results on graphs below; with a slight abuse of notation, we denote the graph $\gamma$ with the same letter as the corresponding diagram.

We say that:

(a) A diagram has a *flat edge* if there is at least one pair $\{(i_1, j_1), (i_2, j_2)\}$ such that $i_1 = i_2$; we write $\gamma \in \Gamma_F(I, J)$ for a diagram with at least a flat edge, and $\gamma \in \Gamma_{\overline{F}}(I, J)$ otherwise. A graph corresponding to a diagram with a flat edge includes an edge $i_k i_k$ which arrives in the same vertex where it started; for these circumstances the term *pseudograph* is preferred by some authors (e.g., [11]).

(b) A diagram $\gamma \in \Gamma_{\overline{F}}(I, J)$ is *connected* if it is not possible to partition the $i$'s into two sets $A, B$ such that there are no edges with $i_1 \in A$ and $i_2 \in B$. We write $\gamma \in \Gamma_C(I, J)$ for connected diagrams, $\gamma \in \Gamma_{\overline{C}}(I, J)$ otherwise. Obviously a diagram is connected if and only if the corresponding graph is connected, in the standard sense.

(c) A diagram $\gamma \in \Gamma_{\overline{F}}(I, J)$ is *paired* if, considering any two sets of edges $\{(i_1, j_1), (i_2, j_2)\}$ and $\{(i_3, j_3), (i_4, j_4)\}$, then $i_1 = i_3$ implies $i_2 = i_4$; in words, the rows are completely coupled two by two. We write $\gamma \in \Gamma_P(I, J)$ for paired diagrams.

It is obvious that for $P > 2$ a paired diagram cannot be connected. Note that if $Q$ is odd, paired diagrams cannot have flat edges, so that the assumption $\gamma \in \Gamma_{\overline{F}}(I, J)$ becomes redundant.

(d) We shall say a diagram has a *k-loop* if there exists a sequence of $k$ edges

$$\{(i_1, j_1), (i_2, j_2)\}, \ldots, \{(i_k, j_k), (i_{k+1}, j_{k+1})\} = (i_1 i_2), \ldots, (i_k i_{k+1})$$

such that $i_1 = i_{k+1}$; we write $\gamma \in \Gamma_{L(k)}(I, J)$ for diagrams with a $k$-loop and no loop of order smaller than $k$.

Note that $\Gamma_F(I, J) = \Gamma_{L(1)}(I, J)$ (a flat edge is a 1-loop); also, we write

$$\Gamma_{CL(k)}(I, J) = \Gamma_C(I, J) \cap \Gamma_{L(k)}(I, J)$$

for connected diagrams with $k$-loops, and $\Gamma_{C\overline{L(k)}}(I, J)$ for connected diagrams with no loops of order $k$ or smaller. For instance, a connected diagram belongs to $\Gamma_{C\overline{L(2)}}(I, J)$ if there are neither flat edges nor two edges $\{(i_1, j_1), (i_2, j_2)\}$ and $\{(i_3, j_3), (i_4, j_4)\}$ such that $i_1 = i_3$ and $i_2 = i_4$; in words, there are no pairs of rows which are connected twice.

A graph is *Hamiltonian* [written $\gamma \in \Gamma_H(I, J)$] [11] if it has a spanning cycle, that is, if there exists a loop which covers all the vertices without touching any of them (other than the first) twice. Two graphs $G_1 = (I_1, E_1)$ and $G_2 = (I_2, E_2)$ are isomorphic if there exists a one-to-one, onto mapping $\phi : I_1 \to I_2$ such that $i_1 i_2 \in E_1 \Longleftrightarrow \phi(i_1)\phi(i_2) \in E_2$.

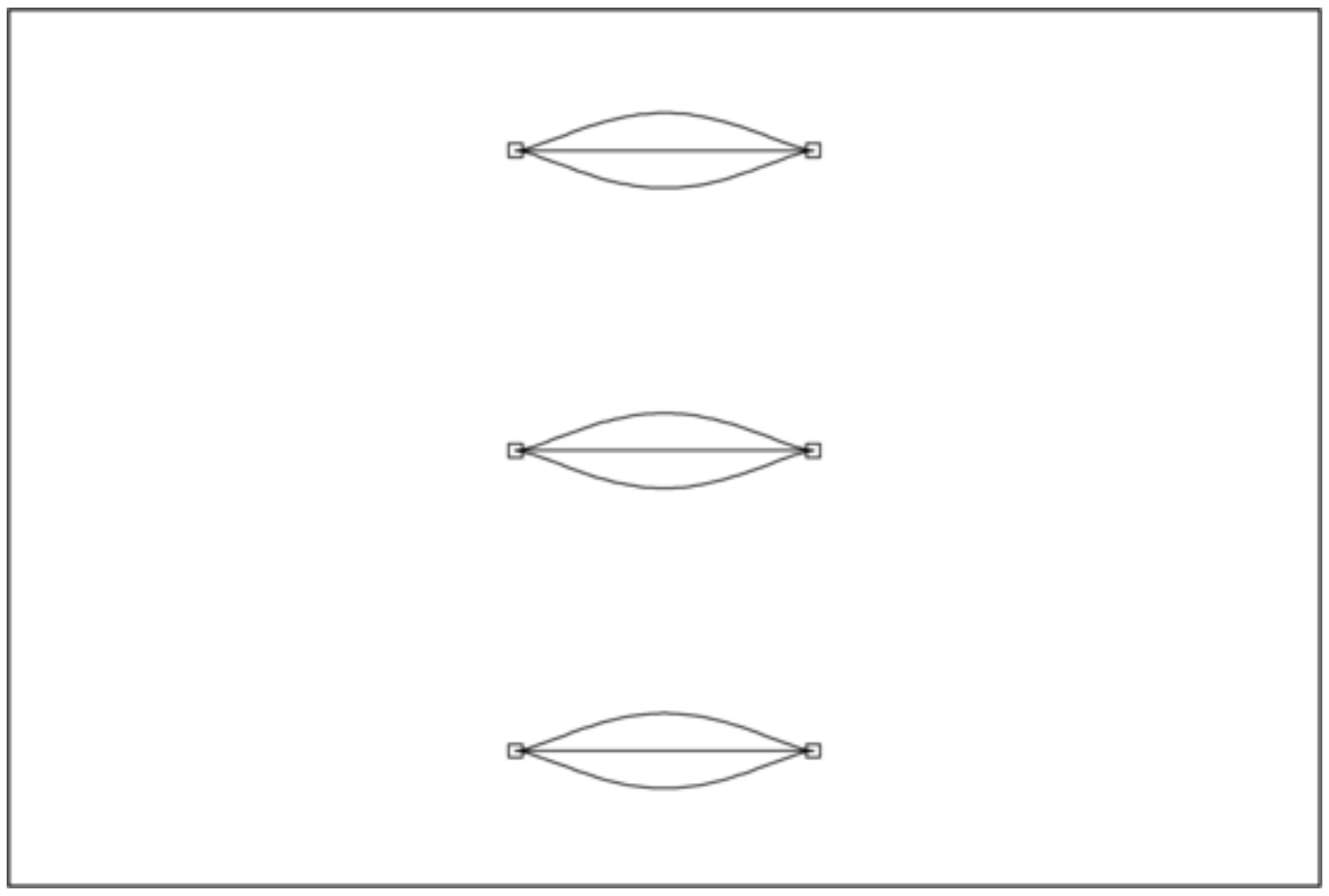

FIG. 1. *A multigraph for* $\gamma \in \Gamma_P(6,3)$.

In many cases either $I$ or $J$ (or both) can be simply taken as the set of the first $p$ or $q$ natural numbers, that is, $I = \{1,\dots,p\}$, $J = \{1,\dots,q\}$. Under such circumstances, when confusion is possible we shall occasionally write $\Gamma(I,q), \Gamma(p,J)$ or $\Gamma(p,q)$ for $\Gamma(I,J)$. Some examples of graphs are drawn in Figures 1–6.

We have the following result.

THEOREM 3.1. *For all* $l_1 \leq l_2 \leq l_3$ *we have*

$$EI^2_{l_1l_2l_3} = \Delta_{l_1l_2l_3}; \tag{10}$$

*moreover, for* $p = 2,3,4$,

$$EI^{2p}_{l_1l_2l_3} = (2p-1)!!\Delta^p_{l_1l_2l_3} + O(l_1^{-1}). \tag{11}$$

PROOF. Result (10) is known in the physics literature; see, for instance, [19]. For (11), we recall that

$$Ea_{l_{j_1}m_{i_1j_1}}a_{l_{j_2}m_{i_2j_2}} = (-1)^{m_{i_1j_1}}C_{l_{j_1}}\delta^{l_{j_2}}_{l_{j_1}}\delta^{-m_{i_2j_2}}_{m_{i_1j_1}}; \tag{12}$$

hence, in view of (9), and because the spherical harmonic coefficients are (complex) Gaussian distributed, the following formula holds, for all $I$:

$$E\left\{\prod_{i\in I}\prod_{j=1}^{3}\frac{a_{l_jm_{ij}}}{\sqrt{C_{l_i}}}\right\} = \sum_{\gamma\in\Gamma(I,3)}\delta(\gamma;l_1,l_2,l_3), \tag{13}$$

where we define

$$\delta(\gamma;l_1,l_2,l_3) = \prod_{\{(i_uj_u),(i'_uj'_u)\}\in\gamma}(-1)^{m_{i_uj_u}}\delta^{-m_{i'_uj'_u}}_{m_{i_uj_u}}\delta^{l_{j'_u}}_{l_{j_u}}. \tag{14}$$

For any diagram $\gamma$, we can also write

$$D[\gamma;l_1,l_2,l_3]=\left\{\prod_{i\in I}\prod_{j=1}^{3}\sum_{m_{ij}=-l_j}^{l_j}\right\}\prod_{i\in I}\begin{pmatrix} l_1 & l_2 & l_3\\ m_{i1} & m_{i2} & m_{i3}\end{pmatrix}\delta(\gamma;l_1,l_2,l_3), \tag{15}$$

where

$$\left\{\prod_{i\in I}\prod_{j=1}^{3}\sum_{m_{ij}=-l_j}^{l_j}\right\}=\sum_{m_{i_1 1}}\cdots\sum_{m_{i_P 3}},\qquad \{i_1,\dots,i_P\}=I;$$

to be more explicit, there are $3\times P$ summations to compute: for instance, when $P=2$ we get

$$\left\{\prod_{i\in\{i_1,i_2\}}\prod_{j=1}^{3}\sum_{m_{ij}=-l_j}^{l_j}\right\}=\sum_{m_{i_1 1}=-l_1}^{l_1}\sum_{m_{i_1 2}=-l_2}^{l_2}\sum_{m_{i_1 3}=-l_3}^{l_3}\sum_{m_{i_2 1}=-l_1}^{l_1}\sum_{m_{i_2 2}=-l_2}^{l_2}\sum_{m_{i_2 3}=-l_3}^{l_3}.$$

Furthermore, we also define

$$\begin{aligned}
&D[\Gamma(I,3);l_1,l_2,l_3]\\
&\quad=\sum_{\gamma\in\Gamma(I,3)}D[\gamma;l_1,l_2,l_3]\\
&\quad=\left\{\prod_{i\in I}\prod_{j=1}^{3}\sum_{m_{ij}=-l_j}^{l_j}\right\}\prod_{i\in I}\begin{pmatrix} l_1 & l_2 & l_3\\ m_{i1} & m_{i2} & m_{i3}\end{pmatrix}\sum_{\gamma\in\Gamma(I,3)}\delta(\gamma;l_1,l_2,l_3);
\end{aligned}$$

in words, $D[\cdot;l_1,l_2,l_3]$ represents the component of the expected value that corresponds to a particular set of diagrams. Notice that

$$\begin{aligned}
EI_{l_1l_2l_3}^{2p}&=\sum_{m_{11}=-l_1}^{l_1}\cdots\sum_{m_{2p,3}=-l_3}^{l_3}E\left\{\prod_{i=1}^{2p}\left[\begin{pmatrix} l_1 & l_2 & l_3\\ m_{i1} & m_{i2} & m_{i3}\end{pmatrix}\prod_{j=1}^{3}\frac{a_{l_jm_{ij}}}{\sqrt{C_{l_j}}}\right]\right\}\\
&=\sum_{m_{11}=-l_1}^{l_1}\cdots\sum_{m_{2p,3}=-l_3}^{l_3}\left\{\prod_{i=1}^{2p}\begin{pmatrix} l_1 & l_2 & l_3\\ m_{i1} & m_{i2} & m_{i3}\end{pmatrix}\right\}E\left\{\prod_{i=1}^{2p}\prod_{j=1}^{3}\frac{a_{l_jm_{ij}}}{\sqrt{C_{l_j}}}\right\}\\
&=\sum_{m_{11}=-l_1}^{l_1}\cdots\sum_{m_{2p,3}=-l_3}^{l_3}\left\{\prod_{i=1}^{2p}\begin{pmatrix} l_1 & l_2 & l_3\\ m_{i1} & m_{i2} & m_{i3}\end{pmatrix}\right\}\sum_{\gamma\in\Gamma(2p,3)}\delta(\gamma;l_1,l_2,l_3)\\
&=D[\Gamma(2p,3);l_1,l_2,l_3].
\end{aligned}$$

Now

$$\begin{aligned}
D[\Gamma(2p,3);l_1,l_2,l_3]&=D[\Gamma_P(2p,3);l_1,l_2,l_3]\\
&\qquad+D[\Gamma(2p,3)\backslash\Gamma_P(2p,3);l_1,l_2,l_3];
\end{aligned}$$

our aim is to show that

$$D[\Gamma_P(2p,3);l_1,l_2,l_3]=(2p-1)!!\Delta^p_{l_1l_2l_3}, \tag{16}$$

$$D[\Gamma(2p,3)\backslash\Gamma_P(2p,3);l_1,l_2,l_3]=O(l_1^{-1}). \tag{17}$$

To establish (16) and (17) we rely on Propositions 3.1 and 3.2, whose proofs are collected in the Appendix. □

The next lemma relates to the "Gaussian" component of the expected value, that is, the diagrams that are paired.

Proposition 3.1. *For any $p\in\mathbb{N}$, and $I$ with cardinality $\#(I)=2p$, we have*

$$D[\Gamma_P(I,3);l_1,l_2,l_3]=(2p-1)!!\Delta^p_{l_1l_2l_3}.$$

The proof that $(17)=O(l_1^{-1})$ requires considerably more work; the next three lemmas refer to diagrams with loops of order 1, 2 and 3, respectively.

Lemma 3.1. *For diagrams with a flat edge, $\gamma\in\Gamma_F(I,3)$, we have*

$$D[\gamma;l_1,l_2,l_3]=0.$$

The next two results show how diagrams belonging to $\Gamma_{CL(2)}(I,3)$, $\Gamma_{CL(3)}(I,3)$ can be "reduced"; namely, they show how the corresponding summands in the expected value can be expressed in terms of smaller-order diagrams. Without loss of generality we can take $\#(I)\geq 4$; indeed the case $\#(I)=2$ has been dealt with in the proof of Theorem 3.1, while we know that odd moments are identically equal to zero. Let $\gamma\in\Gamma_{CL(2)}(I,3)$ be a connected diagram with a 2-loop, and denote by $i_1,i_2$ the rows that are linked by two edges; in other words, the diagram includes both the edges $[(i_1,j_1),(i_2,j_1')]$ and $[(i_1,j_2),(i_2,j_2')]$; in the sequel, $j_k,j_k'$ takes values in $(1,2,3)$, for any integer $k$. Because the diagram is connected, there must exist also edges $[(i_1,j_3),(i_3,j_3')]$ and $[(i_2,j_4),(i_4,j_4')]$, where $i_3,i_4\neq i_1,i_2$. We denote by $\gamma_{R(i_1,i_2)}$ the lower-order diagram which is obtained by deleting $[(i_1,j_1),(i_2,j_1')]$ and $[(i_1,j_2),(i_2,j_2')]$, and substituting $[(i_1,j_3),(i_3,j_3')]$ and $[(i_2,j_4),(i_4,j_4')]$ with $[(i_3,j_3'),(i_4,j_4')]$. In graphical terms, $\gamma_{R(i_1,i_2)}$ is obtained by cutting the two nodes $i_1,i_2$ and merging together the edges that departed from them to reach other vertices; $\gamma_{R(i_1,i_2)}$ can itself belong to $\Gamma_{CL(2)}(I-2,3)$ and the argument can be iterated (Figure 2). We note that these reductions need not be unique in general; any arbitrary choice of a suitable pair of nodes would not affect our argument, however.

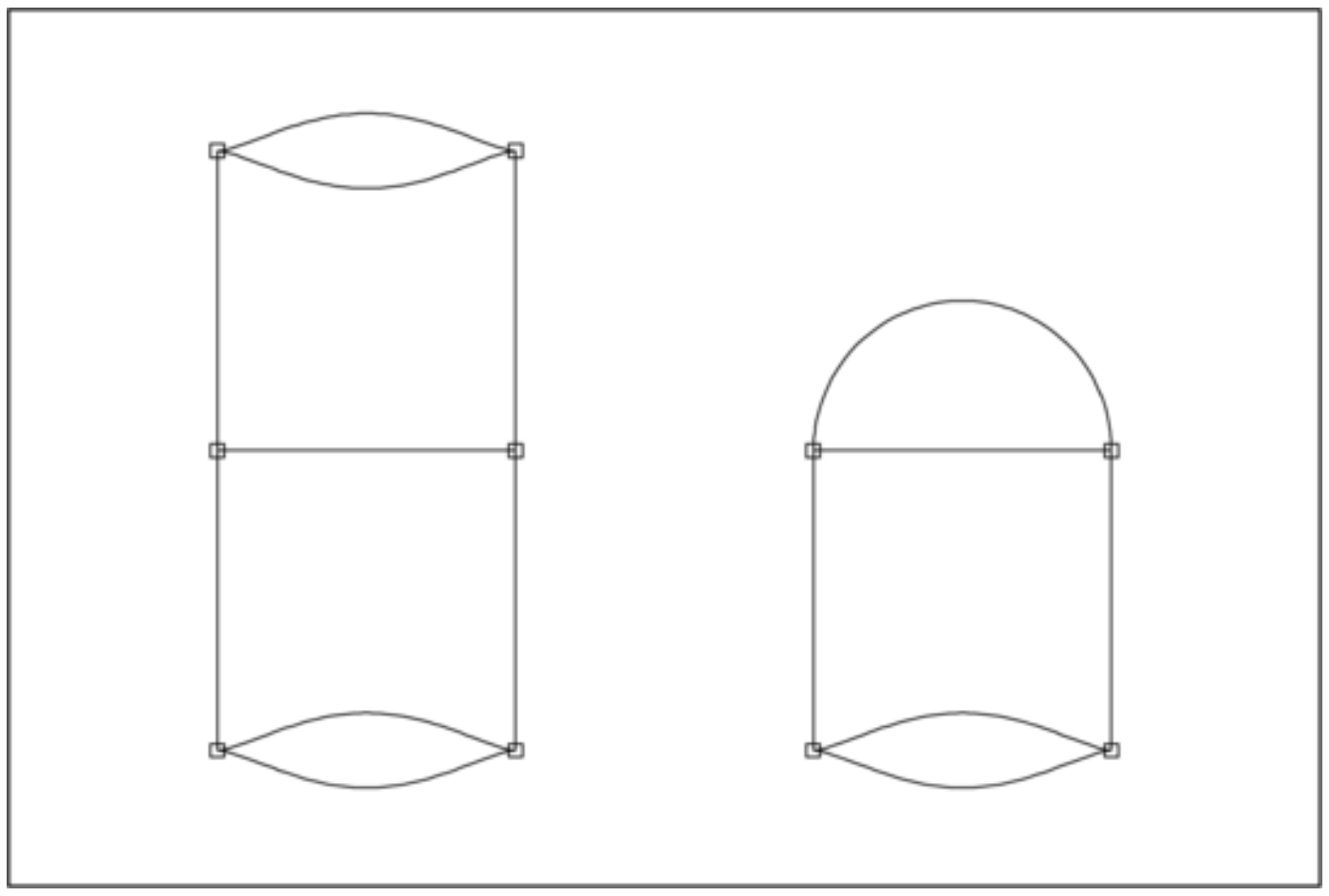

FIG. 2. $\gamma \in \Gamma_{CL(2)}(6,3)$ *and* $\gamma_{R(i_1,i_2)}$.

LEMMA 3.2. *For* $\gamma \in \Gamma_{CL(2)}(I,3)$ *and* $\gamma_{R(i_1,i_2)}$ *as defined before, we have*

$$D[\gamma; l_1, l_2, l_3] = \frac{1}{2l_{j_3}+1} D[\gamma_{R(i_1,i_2)}; l_1, l_2, l_3].$$

By definitions (14) and (15) $D[\gamma; l_1, l_2, l_3]$ can be nonzero only if $l_{j_3} = l_{j'_3} = l_{j_4} = l_{j'_4}$, so there is no notational ambiguity in Lemma 3.2. More explicitly, assuming, for instance, that the edges $[(1,1),(2,1)]$ and $[(1,2),(2,2)]$ are present in $\gamma$, then a factor $(2l_3+1)^{-1}$ will emerge from the reduction.

We now focus on diagrams with a 3-loop. Let $\gamma \in \Gamma_{CL(3)}(I,3)$ be a connected diagram with a 3-loop, and denote by $i_1, i_2, i_3$ the rows that are linked by the loop; in other words, the diagram includes the three edges

$$[(i_1,j_1),(i_2,j_2)], \qquad [(i_2,j_3),(i_3,j_4)], \qquad [(i_3,j_5),(i_1,j_6)].$$

Because the diagram is connected, there must exist also edges

$$[(i_1,j_7),(i_4,j_8)], \qquad [(i_2,j_9),(i_5,j_{10})], \qquad [(i_3,j_{11}),(i_6,j_{12})]$$

where $i_4, i_5, i_6 \neq i_1, i_2, i_3$. We denote by $\gamma_{R(i_1,i_2,i_3)}$ the lower-order diagram which is obtained by replacing $i_2, i_3$ with $i_1$ and then deleting all flat edges. More explicitly, $\gamma_{R(i_1,i_2,i_3)}$ is obtained by deleting

$$[(i_1,j_1),(i_2,j_2)], \qquad [(i_2,j_3),(i_3,j_4)], \qquad [(i_3,j_5),(i_1,j_6)],$$

and substituting

$$[(i_2,j_9),(i_5,j_{10})], \qquad [(i_3,j_{11}),(i_6,j_{12})]$$

with

$$[(i_1,j_9),(i_5,j_{10})], \qquad [(i_1,j_{11}),(i_6,j_{12})].$$

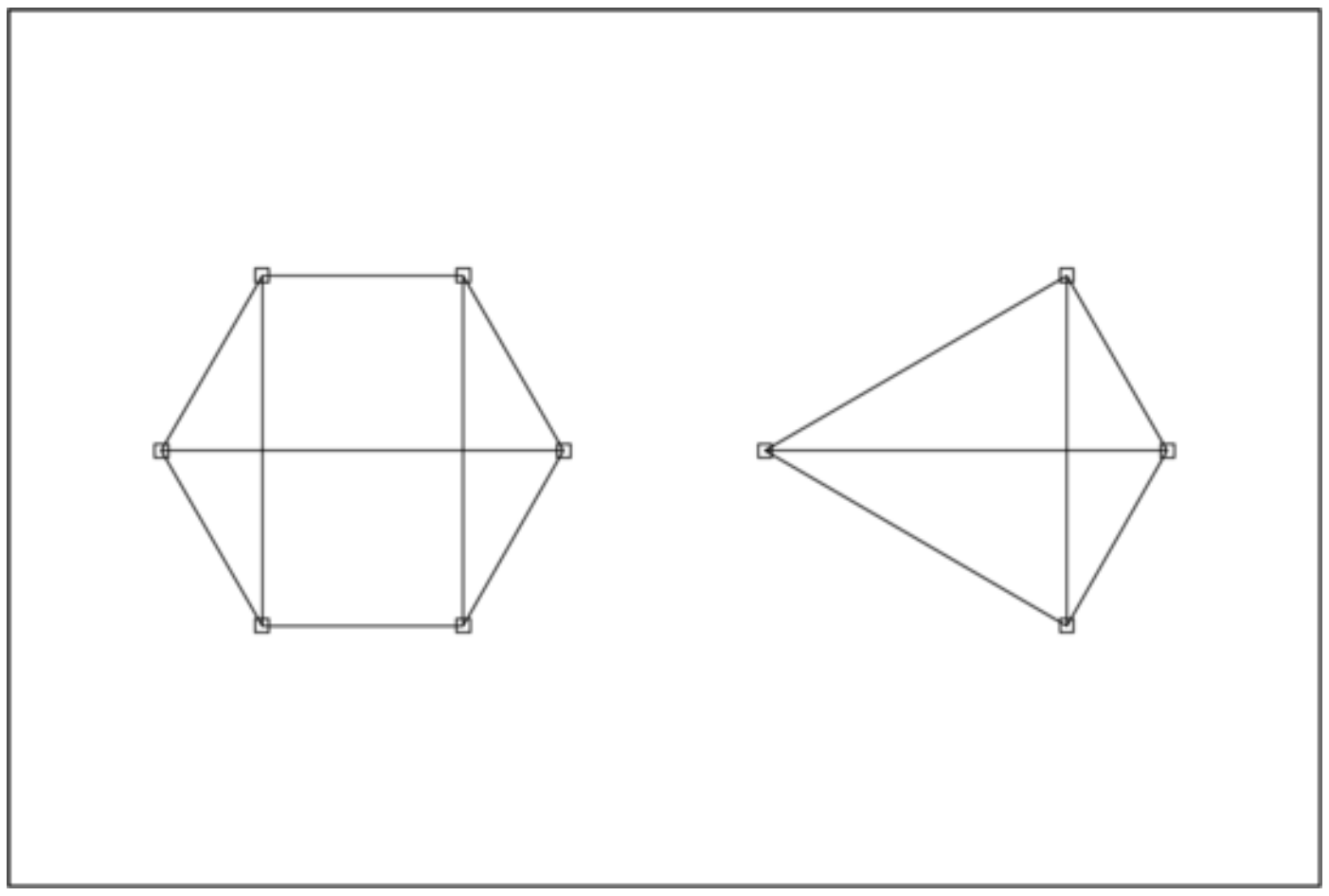

FIG. 3. $\gamma \in \Gamma_{CL(3)}(6,3)$ *and* $\gamma_{R(i_1,i_2,i_3)}$.

In graphical terms, we are merging three nodes into a single one (see Figure 3); again, $\gamma_{R(i_1,i_2,i_3)}$ can belong to $\Gamma_{CL(3)}(I-2,3)$ and the argument can be iterated.

LEMMA 3.3. *For $\gamma \in \Gamma_{CL(3)}(I,3)$ a connected diagram with a 3-loop and $\gamma_{R(i_1,i_2,i_3)}$ as defined before, we have*

$$D[\gamma; l_1, l_2, l_3] = \begin{Bmatrix} l_1 & l_2 & l_3 \\ l_1 & l_2 & l_3 \end{Bmatrix} D[\gamma_{R(i_1,i_2,i_3)}; l_1, l_2, l_3],$$

*where on the left-hand side we have used Wigner's $6j$ coefficient, defined in the Appendix; hence*

$$D[\gamma; l_1, l_2, l_3] = O(l_3^{-1} D[\gamma_{R(i_1,i_2,i_3)}; l_1, l_2, l_3]).$$

The proofs of Lemmas 3.1, 3.2 and 3.3 can be obtained as an application of the graphical method described, for instance, in ([34], Chapters 11 and 12), and they are hence omitted for brevity's sake.

The next proposition exploits the previous results to provide a bound (17) on the "non-Gaussian" part of the higher-order moments of the angular bispectrum.

PROPOSITION 3.2. *For all $I$ such that $\#(I) = 4, 6$ or $8$ and $l_1 \leq l_2 \leq l_3$, we have*

$$D[\Gamma(I,3) \backslash \Gamma_P(I,3); l_1, l_2, l_3] = O(l_1^{-1}).$$

REMARK 3.1. A careful inspection of the proofs of Propositions 3.1 and 3.2 reveals that for $l_1 < l_2 < l_3$ we obtain the special case

$$
\begin{aligned}
EI^4_{l_1l_2l_3} &= 3 + \frac{6}{2l_1+1} + \frac{6}{2l_2+1} + \frac{6}{2l_3+1} + 6\begin{Bmatrix} l_1 & l_2 & l_3 \\ l_1 & l_2 & l_3 \end{Bmatrix} \\
&= \frac{6l_1+9}{2l_1+1} + O(l_2^{-1}),
\end{aligned}
\tag{18}
$$

a result that we shall exploit for statistical applications.

REMARK 3.2. It is natural to conjecture that a result analogous to Proposition 3.2 will hold for all sets $\#(I) = 2p$, $p \in N$. This conjecture is not simple to explore, however, because the proofs in the Appendix require some analytic properties of the summations of (products of) Wigner's $3j$ coefficients, which have not been extended (to the best of our knowledge) to products of arbitrary order.

**4. Unknown angular power spectrum.** In this section we focus on the more realistic case where the angular power spectrum is unknown and estimated from the data; so we consider $\widehat{I}_{l_1l_2l_3}$ rather than $I_{l_1l_2l_3}$. As before, under Gaussianity of the underlying field $T(\theta,\varphi)$

$$
E\widehat{I}^{2p-1}_{l_1l_2l_3} = 0, \qquad p = 1, 2, \ldots,
$$

by a simple symmetry argument. Now note that

$$
\begin{aligned}
&\left(\frac{|a_{l0}|^2}{\widehat{C}_l}, \frac{2|a_{l1}|^2}{\widehat{C}_l}, \ldots, \frac{2|a_{ll}|^2}{\widehat{C}_l}\right) \\
&\quad = (2l+1)\left(\frac{|a_{l0}|^2}{|a_{l0}|^2 + \sum_{m=1}^{l} 2|a_{lm}|^2}, \frac{2|a_{l1}|^2}{|a_{l0}|^2 + \sum_{m=1}^{l} 2|a_{lm}|^2}, \ldots, \right. \\
&\qquad\qquad \left. \frac{2|a_{ll}|^2}{|a_{l0}|^2 + \sum_{m=1}^{l} 2|a_{lm}|^2}\right) \\
&\quad \overset{\text{def}}{=} (2l+1)(\xi_{l0}, \ldots, \xi_{ll}) \overset{d}{=} (2l+1)\,\mathrm{Dir}\left(\frac{1}{2}, 1, \ldots, 1\right);
\end{aligned}
$$

here $\overset{d}{=}$ denotes equality in distribution and $\mathrm{Dir}(\theta_0, \ldots, \theta_p)$ a Dirichlet distribution with parameters $(\theta_0, \ldots, \theta_p)$. Define

$$
u_{lm} = \frac{a_{lm}}{\sqrt{C_l}}, \qquad \widehat{u}_{lm} = \frac{a_{lm}}{\sqrt{\widehat{C}_l}}, \qquad m = 0, 1, \ldots, l.
\tag{19}
$$

We have the following simple result.

THEOREM 4.1. *Let $l$ and $p$ be positive integers, and define*

$$g(l;p)=\prod_{k=1}^{p}\left\{\frac{2l+1}{2l+2k-1}\right\}.$$

*Now for $u$ and $\hat{u}$ defined by* (19), *we have*

$$\begin{aligned}
&E\Bigl\{\underbrace{\hat{u}_{l0}\dots\hat{u}_{l0}}_{q_0\text{ times}}\underbrace{\hat{u}_{l1}\dots\hat{u}_{l1}}_{q_1\text{ times}}\underbrace{\hat{u}^*_{l1}\dots\hat{u}^*_{l1}}_{q'_1\text{ times}}\cdots\underbrace{\hat{u}_{lk}\dots\hat{u}_{lk}}_{q_k\text{ times}}\underbrace{\hat{u}^*_{lk}\dots\hat{u}^*_{lk}}_{q'_k\text{ times}}\Bigr\}\\
&\quad=E\Bigl\{\underbrace{u_{l0}\dots u_{l0}}_{q_0\text{ times}}\underbrace{u_{l1}\dots u_{l1}}_{q_1\text{ times}}\underbrace{u^*_{l1}\dots u^*_{l1}}_{q'_1\text{ times}}\cdots\underbrace{u_{lk}\dots u_{lk}}_{q_k\text{ times}}\underbrace{u^*_{lk}\dots u^*_{lk}}_{q'_k\text{ times}}\Bigr\}\\
&\qquad\times g(l;q_0+q_1+\cdots+q'_k).
\end{aligned}$$

PROOF. By symmetry arguments, it is easy to see that both sides are zero unless $q_0=2p_0$ (say) is even and $q_i=q'_i=p_i$ (say), for $i=1,\dots,k$. The $u_{lm}$ are independent over different $m$'s, and thus we have

$$\begin{aligned}
&E\Bigl\{\underbrace{u_{l0}\dots u_{l0}}_{2p_0\text{ times}}\underbrace{u_{l1}\dots u_{l1}}_{q_1\text{ times}}\underbrace{u^*_{l1}\dots u^*_{l1}}_{q'_1\text{ times}}\cdots\underbrace{u_{lk}\dots u_{lk}}_{q_k\text{ times}}\underbrace{u^*_{lk}\dots u^*_{lk}}_{q'_k\text{ times}}\Bigr\}\\
&\quad=Eu_{l0}^{2p_0}Eu_{l1}^{p_1}(u^*_{l1})^{p_1}\cdots Eu_{lk}^{p_k}(u^*_{lk})^{p_k}\\
&\quad=\begin{cases}(2p_0-1)!!\prod_{i=1}^{k}p_i!, & \text{for } p_0>0,\\ \prod_{i=1}^{k}p_i!g(l;p), & \text{for } p_0=0,\end{cases}
\end{aligned}$$

because

$$u_{l0}\stackrel{d}{=}N(0,1)\quad\text{and}\quad u_{lm}u^*_{lm}=|u_{lm}|^2\stackrel{d}{=}\exp(1).$$

Now write $p=p_0+\cdots+p_k$, and note that ([18], page 233)

$$\begin{aligned}
&E\{\hat{u}_{l0}^{2p_0}\hat{u}_{l1}^{p_1}(\hat{u}^*_{l1})^{p_1}\cdots\hat{u}_{lk}^{p_k}(\hat{u}^*_{lk})^{p_k}\}\\
&\quad=\frac{(2l+1)^p}{2^{p_1+\cdots+p_k}}E\xi_{l0}^{p_0}\cdots\xi_{lk}^{p_k}\\
&\quad=\frac{(2l+1)^p}{2^{p_1+\cdots+p_k}}\frac{\Gamma(l+1/2)}{\Gamma(l+p+1/2)}\frac{\Gamma(p_0+1/2)\Gamma(p_2+1)\cdots\Gamma(p_k+1)}{\Gamma(1/2)}\\
&\quad=\frac{(2p_0-1)!!p_1!\times\cdots\times p_k!}{(2l+1)\times\cdots\times(2l+2p-1)}(2l+1)^p
\end{aligned}$$

$$= \begin{cases} (2p_0-1)!!\displaystyle\prod_{i=1}^{k} p_i! g(l;p), & \text{for } p_0>0, \\ \displaystyle\prod_{i=1}^{k} p_i! g(l;p), & \text{for } p_0=0, \end{cases}$$

as claimed. □

Some special cases are as follows:

$$E\frac{|a_{l0}|^2}{\widehat{C}_l} = E\frac{|a_{lm}|^2}{\widehat{C}_l} = 1, \qquad m=\pm 1,\dots,\pm l,$$

$$E\left\{\frac{|a_{l0}|^2}{\widehat{C}_l}\right\}^p = \frac{(2p-1)!!}{(2l+1)\times\cdots\times(2l+2p-1)}(2l+1)^p,$$

$$E\left\{\frac{|a_{lm}|^2}{\widehat{C}_l}\right\}^p = \frac{(2l+1)^p}{2^p}\frac{\Gamma(l+1/2)}{\Gamma(l+p+1/2)}\frac{\Gamma(1/2)\Gamma(p+1)}{\Gamma(1/2)}$$
$$= \frac{p!}{(2l+1)\times\cdots\times(2l+2p-1)}(2l+1)^p,$$

and for $p=p_1+p_2$, $p_1,p_2>0$

$$E\left\{\left(\frac{|a_{l0}|^2}{\widehat{C}_l}\right)^{p_1}\left(\frac{|a_{l1}|^2}{\widehat{C}_l}\right)^{p_2}\right\} = \frac{(2l+1)^p}{2^p}\frac{\Gamma(l+1/2)}{\Gamma(l+p+1/2)}\frac{\Gamma(p_1+1/2)\Gamma(p_2+1)}{\Gamma(1/2)}$$
$$= \frac{(2p_1-1)!!p_2!}{(2l+1)\times\cdots\times(2l+2p-1)}(2l+1)^p,$$

$$E\left\{\left(\frac{|a_{l1}|^2}{\widehat{C}_l}\right)^{p_1}\left(\frac{|a_{l2}|^2}{\widehat{C}_l}\right)^{p_2}\right\} = \frac{(2l+1)^p}{2^p}\frac{\Gamma(l+1/2)}{\Gamma(l+p+1/2)}$$
$$\times\frac{\Gamma(1/2)\Gamma(p_1+1)\Gamma(p_2+1)}{\Gamma(1/2)}$$
$$= \frac{p_1!p_2!}{(2l+1)\times\cdots\times(2l+2p-1)}(2l+1)^p.$$

By Theorem 4.1, it is possible to establish a simple relationship between the normalized bispectrum with known or unknown angular power spectrum. More precisely, it is immediate to see that for $l_1<l_2<l_3$

$$E\widehat{I}^{2p}_{l_1l_2l_3} = EI^{2p}_{l_1l_2l_3}\prod_{i=1}^{3} g(l_i;p),$$

that is, for instance,

$$E\hat{I}^2_{l_1l_2l_3} = EI^2_{l_1l_2l_3} = 1 \tag{20}$$

and

$$E\hat{I}^4_{l_1l_2l_3} = EI^4_{l_1l_2l_3}\left(1 - \frac{2}{2l_1+3}\right)\left(1 - \frac{2}{2l_2+3}\right)\left(1 - \frac{2}{2l_3+3}\right). \tag{21}$$

Also, for $l_1 = l_2 < l_3$ and $l_1 < l_2 = l_3$

$$E\hat{I}^{2p}_{l_1l_1l_3} = EI^{2p}_{l_1l_1l_3} g(l_1; 2p) g(l_3; p), \qquad E\hat{I}^{2p}_{l_1l_3l_2} = EI^{2p}_{l_1l_3l_3} g(l_1; p) g(l_3; 2p);$$

finally, for $l_1 = l_2 = l_3 = l$

$$E\hat{I}^{2p}_{lll} = EI^{2p}_{lll} g(l; 3p),$$

so that, for instance,

$$E\hat{I}^2_{lll} = 6\left(1 - \frac{2}{2l+3}\right)\left(1 - \frac{4}{2l+5}\right). \tag{22}$$

It is interesting to note that

$$E\hat{I}^{2p}_{l_1l_2l_3} \le EI^{2p}_{l_1l_2l_3} \quad \text{and} \quad \lim_{l_1\to\infty} \frac{E\hat{I}^{2p}_{l_1l_2l_3}}{EI^{2p}_{l_1l_2l_3}} = 1 \tag{23}$$

for all choices of $(l_1, l_2, l_3)$ for which the bispectrum is well defined.

**5. Some statistical applications.** In this section we exploit the previous results to derive the asymptotic convergence of some functionals of the bispectrum array. Because the expected bispectrum is identically zero for Gaussian fields, these functionals arise as natural candidates in the development of statistical tests for non-Gaussianity.

More precisely, assume that the resolution of the experiment is such that it yields a maximum observable multipole equal to $L$. It is, in practice, infeasible to take into account all available bispectrum ordinates for the implementation of a statistical procedure: indeed, for current experiments these ordinates are of the order of $L^3 \sim 10^8/10^9$, and the evaluation of all these statistics is beyond the power of the fastest supercomputers for the near future. It is therefore mandatory to consider only a subset of bispectrum ordinates for the test. There are, of course, several possible choices of configurations. We shall restrict our attention to two of them; precisely, for finite integers $l_0 \ge 2$, $K \ge 0$ we shall consider the processes

$$J_{1L;l_0,K}(r) = \frac{1}{\sqrt{L/2}} \sum_{l \,\text{even}, l=l_0+K}^{[Lr]} \left\{ \frac{1}{\sqrt{K+1}} \sum_{u=0}^{K} \frac{\hat{I}_{l-u,l,l+u}}{\sqrt{\Delta_{l-u,l,l+u}}} \right\}, \tag{24}$$

$$J_{2L;l_0,K}(r) = \frac{1}{\sqrt{L/2}} \sum_{l \,\text{even}, l=l_0+K}^{[Lr]} \left\{ \frac{1}{\sqrt{K+1}} \sum_{u=0}^{K} \frac{(\hat{I}^2_{l-u,l,l+u} - \Delta_{l-u,l,l+u})}{\sqrt{2}\Delta_{l-u,l,l+u}} \right\} \tag{25}$$

and

$$J_{3L;l_0,K}(r) = \frac{1}{\sqrt{L}} \sum_{l=l_0+K+1}^{[Lr]-l_0-K} \left\{ \frac{1}{\sqrt{K+1}} \sum_{u=0}^{K} \hat{I}_{l_0+u,l,l+l_0+u} \right\}, \tag{26}$$

$$J_{4L;l_0,K}(r) = \frac{1}{\sqrt{L}} \sum_{l=l_0+K+1}^{[Lr]-l_0-K} \left\{ \frac{1}{\sqrt{K+1}} \sum_{u=0}^{K} \frac{\hat{I}^2_{l_0+u,l,l+l_0+u} - 1}{\sqrt{2}} \right\}, \tag{27}$$

where $[\cdot]$ denotes the integer part of a real number; $0 \le r \le 1$ and $l_0$ is an (arbitrary but fixed) value which can be taken equal to 2, for instance, for cosmological applications (remember the dipole $l_0 = 1$ is usually discarded from CMB data, as it is associated with kinematic effects mainly due to the motion of the Milky Way and the local group of galaxies). As usual, the sums are taken to be equal to zero when the index set is empty. $K$ is a fixed pooling parameter: for $K = 0$ we obtain the special cases

$$\begin{aligned} J_{1L;l_0}(r) &= \frac{1}{\sqrt{L/2}} \sum_{l \,\text{even}, l \ge l_0}^{[Lr]} \frac{\hat{I}_{lll}}{\sqrt{6}}, \\ J_{2L;l_0}(r) &= \frac{1}{\sqrt{L/2}} \sum_{l \,\text{even}, l \ge l_0}^{[Lr]} \frac{(\hat{I}^2_{lll} - 6)}{6\sqrt{2}} \end{aligned} \tag{28}$$

and

$$\begin{aligned} J_{3L;l_0}(r) &= \frac{1}{\sqrt{L}} \sum_{l=l_0+1}^{[Lr]-l_0} \hat{I}_{l_0,l,l+l_0}, \\ J_{4L;l_0}(r) &= \frac{1}{\sqrt{L}} \sum_{l=l_0+1}^{[Lr]-l_0} \frac{\hat{I}^2_{l_0,l,l+l_0} - 1}{\sqrt{2}}. \end{aligned} \tag{29}$$

The normalizing factors are chosen to ensure an asymptotic unit variance for all summands, by means of Theorem 3.1 and (18), (21) and (22). For instance,

$$\begin{aligned} &\mathrm{Var}\{\hat{I}^2_{l_0,l,l+l_0} - 1\} \\ &\quad = \mathrm{Var}\{\hat{I}^2_{l_0,l,l+l_0}\} = E\hat{I}^4_{l_0,l,l+l_0} - \{E\hat{I}^2_{l_0,l,l+l_0}\}^2 \\ &\quad = EI^4_{l_0,l,l+l_0} \left(\frac{2l_0+1}{2l_0+3}\right)\left(\frac{2l+1}{2l+3}\right)\left(\frac{2l+2l_0+1}{2l+2l_0+3}\right) - 1 \\ &\quad = \left(\frac{6l_0+9}{2l_0+1} + O(l^{-1})\right)\left\{\frac{2l_0+1}{2l_0+3}(1+O(l^{-1}))\right\} - 1 \\ &\quad = 2 + O(l^{-1}). \end{aligned} \tag{30}$$

The two pairs of processes $J_{1L;l_0,K}(r), J_{2L;l_0,K}(r)$ and $J_{3L;l_0,K}(r), J_{4L;l_0,K}(r)$ can be viewed as sorts of boundary cases for the possible configurations of multipoles. Indeed, although none of them has so far been considered in the literature (to the best of our knowledge), it seems natural to view $J_{1L;l_0,K}(r), J_{2L;l_0,K}(r)$ as proposals very close to much of what has been done so far in CMB data analysis. More precisely, it has become very common to restrict attention to multipoles close to or on the "main diagonal" $l_1 = l_2 = l_3 = l$, under the unproved conjecture that the greatest part of the non-Gaussian signal should concentrate in that area. This same rationale motivates $J_{1L;l_0,K}(r), J_{2L;l_0,K}(r)$; we shall show below, though, how such a choice can be very far from optimal in relevant cases. On the other hand, $J_{3L;l_0,K}(r), J_{4L;l_0,K}(r)$ rely on a sort of opposite strategy, that is, for a fixed $l_0$ we aim at maximizing the distance among multipoles, albeit preserving the triangle conditions $l_i \leq l_j + l_k$. There are several alternative procedures one may wish to consider, but those we mentioned lend themselves to a simple analysis, while highlighting some quite unexpected features of asymptotics for fixed-radius fields.

THEOREM 5.1. *As $L \to \infty$, for any fixed integers $l_0 > 0$, $K \geq 0$,*

$$J_{1L;l_0,K}(r), J_{2L;l_0,K}(r), J_{3L;l_0,K}(r), J_{4L;l_0,K}(r) \Rightarrow W(r), \qquad 0 \leq r \leq 1, \tag{31}$$

*where $\Rightarrow$ denotes weak convergence in the Skorohod space $D[0,1]$ and $W(r)$ denotes standard Brownian motion.*

PROOF. The proofs for the processes $J_{aL;l_0,K}(r)$, $a = 1, \dots, 4$, are very similar; we give the details only for the most difficult case, namely $J_{4L;l_0,K}(r)$. Here the proof is made harder by the complicated structure of dependence; note indeed that the set of random coefficients $\{a_{lm} : l = l_0, \dots, l_0 + K, m = -l, \dots, l\}$ belongs to each summand in (29). Denote by $\Im_l$ the filtration generated by the triangular array $\{a_{l,-l}, \dots, a_{l,l}\}$, $l = 1, 2, \dots$, and define

$$X_{l,L} = \frac{1}{\sqrt{K+1}} \sum_{u=0}^{K} \frac{\hat{I}^2_{l_0+u,l,l+l_0+u} - 1}{\sqrt{2L}}, \qquad l = l_0 + K + 1, l_0 + K + 2, \dots,$$

that is,

$$J_{4L;l_0,K}(r) = \sum_{l=l_0+K+1}^{[Lr]-l_0-K} X_{l,L}.$$

Now we note first that

$$E\{X_{l,L} | \Im_{l-m}\} = \frac{1}{\sqrt{K+1}} \sum_{u=0}^{K} \frac{E\{\hat{I}^2_{l_0+u,l,l+l_0+u} | \Im_{l-m}\} - 1}{\sqrt{2L}} = 0, \qquad m \geq 1, \tag{32}$$

because for all $0 < l_1 < l < l_2 < l_3$,

$$
\begin{aligned}
E\{\hat{I}^2_{l_1 l_2 l_3} | \Im_l\} &= \sum_{m_1,m_2,m_3} \sum_{m_1',m_2',m_3'} \begin{pmatrix} l_1 & l_2 & l_3 \\ m_1 & m_2 & m_3 \end{pmatrix} \begin{pmatrix} l_1 & l_2 & l_3 \\ m_1' & m_2' & m_3' \end{pmatrix} \\
&\qquad\qquad \times E\left\{ \frac{a_{l_1 m_1} a_{l_1 m_1'} a_{l_2 m_2} a_{l_2 m_2'} a_{l_3 m_3} a_{l_3 m_3'}}{\widehat{C}_{l_1} \widehat{C}_{l_2} \widehat{C}_{l_3}} \Big| \Im_l \right\} \\
&= \sum_{m_1,m_2,m_3} \sum_{m_1',m_2',m_3'} \begin{pmatrix} l_1 & l_2 & l_3 \\ m_1 & m_2 & m_3 \end{pmatrix} \begin{pmatrix} l_1 & l_2 & l_3 \\ m_1' & m_2' & m_3' \end{pmatrix} \\
&\qquad\qquad \times \frac{a_{l_1 m_1} a_{l_1 m_1'}}{\widehat{C}_1} E\left\{ \frac{a_{l_2 m_2} a_{l_2 m_2'} a_{l_3 m_3} a_{l_3 m_3'}}{\widehat{C}_{l_2} \widehat{C}_{l_3}} \right\} \\
&= \sum_{m_1,m_1'} \sum_{m_2,m_3} \begin{pmatrix} l_1 & l_2 & l_3 \\ m_1 & m_2 & m_3 \end{pmatrix} \begin{pmatrix} l_1 & l_2 & l_3 \\ m_1' & -m_2 & -m_3 \end{pmatrix} \frac{a_{l_1 m_1} a_{l_1 m_1'}}{\widehat{C}_{l_1}} \\
&= \sum_{m_1} \frac{1}{2l_1 + 1} \frac{|a_{l_1 m_1}|^2}{\widehat{C}_{l_1}} = 1.
\end{aligned}
$$

Equation (32) does not imply that the triangular array $\{X_{l,L}\}_{l=2,3,\dots}$ obeys a martingale difference property, because the sequence $X_{l,L}$ is not adapted to the filtration $\Im_l$. However, (32) proves that the pair sequences $\{X_{l,L}, \Im_l\}_{l=2,3,\dots}$ do satisfy a mixingale property [7, 25], that is,

$$
[E(E\{X_{l,L} | \Im_{l-m}\})^2]^{1/2} \le c_1 \frac{m^{-\phi}}{\sqrt{L}} \qquad \text{for } m \ge 1, \tag{33}
$$

$$
[E(X_{l,L} - E\{X_{l,L} | \Im_{l+m}\})^2]^{1/2} \le c_2 \frac{m^{-\phi}}{\sqrt{L}} \qquad \text{for } m \ge 1, \tag{34}
$$

$$
\qquad \text{for some } c_1, c_2, \phi > 0. \tag{35}
$$

Actually the left-hand sides of (33) and (34) are identically zero for $m$ larger than $l_0 + K$, so that, for suitable choices of the constants $c_1, c_2$, the bounds on the right-hand sides hold for an arbitrarily large $\phi$. Note that $\{X^2_{l,L}, l = 1, 2, \dots, L, L = 1, 2, \dots\}$ is a uniformly integrable set, because $\hat{I}^2_{l_1 l_2 l_3}$ has finite fourth-order moments which are uniformly bounded [Theorem 3.1 and (23)]. Also, it is readily seen that

$$
\sup_{0 \le r_1 < r_2 \le 1} \limsup_{L \to \infty} \frac{\sum_{l=[Lr_1]}^{[Lr_2]} E X^2_{l,L}}{r_2 - r_1} < \infty, \qquad \lim_{L \to \infty} \max_{l=1,\dots,L} E X^2_{l,L} = 0.
$$

We have thus established conditions (2.2), (2.3) and (2.5) in [25] for the functional central limit theorem to hold. To complete the proof, we only

need to show that

$$\lim_{L\to\infty} E\left| E\left\{ \left( \sum_{l=[Ls]}^{[Lt]} X_{l,L} \right)^2 \Big| \Im_{[Lr]} \right\} - (t-s) \right| = 0 \qquad \text{for any } r < s < t.$$

For notational simplicity and without loss of generality, we consider the special case $K = 0$. We have

$$\begin{aligned} E&\left\{ \left( \sum_{l=[Ls]}^{[Lt]} X_{l,L} \right)^2 \Big| \Im_{[Lr]} \right\} \\ &= \sum_{l=[Ls]}^{[Lt]} E\{X_{l,L}^2 | \Im_{[Lr]}\} + 2 \sum_{l=[Ls]}^{[Lt]} \sum_{l'=l+1}^{[Lt]} E\{X_{l,L} X_{l',L} | \Im_{[Lr]}\}. \end{aligned}$$

Now for $l, l' > [Lr]$ we have

$$\begin{aligned} E&\{X_{l,L} X_{l',L} | \Im_{[Lr]}\} \\ &= \frac{1}{2L} E\{(\hat{I}^2_{l_0,l,l+l_0} - 1)(\hat{I}^2_{l_0,l',l'+l_0} - 1) | \Im_{[Lr]}\} \\ &= \frac{1}{2L} [E\{\hat{I}^2_{l_0,l,l+l_0} \widehat{I}^2_{l_0,l',l'+l_0} | \Im_{[Lr]}\} \\ &\qquad - E\{\hat{I}^2_{l_0,l,l+l_0} | \Im_{[Lr]}\} - E\{\hat{I}^2_{l_0,l',l'+l_0} | \Im_{[Lr]}\} + 1] \\ &= \frac{1}{2L} [E\{\hat{I}^2_{l_0,l,l+l_0} \widehat{I}^2_{l_0,l',l'+l_0} | \Im_{[Lr]}\} - 1]. \end{aligned} \tag{36}$$

There are now two possible cases, namely $l + l_0 = l'$ and $l + l_0 \neq l'$; in the latter

$$\begin{aligned} E&\{\hat{I}^2_{l_0,l,l+l_0} \hat{I}^2_{l_0,l',l'+l_0} | \Im_{[Lr]}\} \\ &= \sum_{m_{01},m_{02},m_{03},m_{04}} \sum_{m_{11},m_{12},m_{21},m_{22}} \begin{pmatrix} l_0 & l & l+l_0 \\ m_{01} & m_{11} & m_{21} \end{pmatrix} \\ &\quad \times \begin{pmatrix} l_0 & l & l+l_0 \\ m_{02} & -m_{11} & -m_{21} \end{pmatrix} \\ &\quad \times \begin{pmatrix} l_0 & l' & l'+l_0 \\ m_{03} & m_{12} & m_{22} \end{pmatrix} \begin{pmatrix} l_0 & l' & l'+l_0 \\ m_{04} & -m_{12} & -m_{22} \end{pmatrix} \\ &\quad \times \frac{a_{l_0 m_{01}} a_{l_0 m_{02}} a_{l_0 m_{03}} a_{l_0 m_{04}}}{\widehat{C}^2_{l_0}} \\ &= \sum_{m_{01},m_{02},m_{03},m_{04}} \frac{\delta_{m_{01}}^{-m_{02}} \delta_{m_{03}}^{-m_{04}}}{(2l_0+1)^2} \frac{a_{l_0 m_{01}} a_{l_0 m_{02}} a_{l_0 m_{03}} a_{l_0 m_{04}}}{\widehat{C}^2_{l_0}} = 1, \end{aligned} \tag{37}$$

the second to last step following from (49) in the Appendix. Otherwise, for $l+l_0=l'$,

$$
\begin{aligned}
&E\{\hat{I}^2_{l_0,l,l+l_0}\hat{I}^2_{l_0,l',l'+l_0}|\Im_{[Lr]}\}\\
&\quad=\sum_{m_{01},m_{02},m_{03},m_{04}}\sum_{m_{11},m_{12},m_{21},m_{22}}\begin{pmatrix} l_0 & l & l' \\ m_{01} & m_{11} & m_{21}\end{pmatrix}\\
&\qquad\times\begin{pmatrix} l_0 & l & l' \\ m_{02} & -m_{11} & -m_{21}\end{pmatrix}\begin{pmatrix} l_0 & l' & l'+l_0 \\ m_{03} & m_{12} & m_{22}\end{pmatrix}\\
&\qquad\times\begin{pmatrix} l_0 & l' & l'+l_0 \\ m_{04} & -m_{12} & -m_{22}\end{pmatrix}\frac{a_{l_0m_{01}}a_{l_0m_{02}}a_{l_0m_{03}}a_{l_0m_{04}}}{\widehat{C}^2_{l_0}}\\
&\qquad+2\sum_{m_{01},m_{02},m_{03},m_{04}}\sum_{m_{11},m_{12},m_{21},m_{22}}\begin{pmatrix} l_0 & l & l' \\ m_{01} & m_{11} & m_{21}\end{pmatrix}\\
&\qquad\quad\times\begin{pmatrix} l_0 & l & l' \\ m_{02} & m_{12} & -m_{21}\end{pmatrix}\begin{pmatrix} l_0 & l' & l'+l_0 \\ m_{03} & -m_{11} & m_{22}\end{pmatrix}\\
&\qquad\quad\times\begin{pmatrix} l_0 & l' & l'+l_0 \\ m_{04} & -m_{12} & -m_{22}\end{pmatrix}\frac{a_{l_0m_{01}}a_{l_0m_{02}}a_{l_0m_{03}}a_{l_0m_{04}}}{\widehat{C}^2_{l_0}}\\
&\quad=1+2B,
\end{aligned}
\tag{38}
$$

where

$$
\begin{aligned}
B&=\sum_{m_{01},m_{02},m_{03},m_{04}}\sum_{m_{11},m_{12},m_{21},m_{22}}\begin{pmatrix} l_0 & l & l' \\ m_{01} & m_{11} & m_{21}\end{pmatrix}\begin{pmatrix} l_0 & l & l' \\ m_{02} & m_{12} & -m_{21}\end{pmatrix}\\
&\quad\times\begin{pmatrix} l_0 & l' & l'+l_0 \\ m_{03} & -m_{11} & m_{22}\end{pmatrix}\begin{pmatrix} l_0 & l' & l'+l_0 \\ m_{04} & -m_{12} & -m_{22}\end{pmatrix}\frac{a_{l_0m_{01}}a_{l_0m_{02}}a_{l_0m_{03}}a_{l_0m_{04}}}{\widehat{C}^2_{l_0}}\\
&=\sum_{s=0}^{\infty}\sum_{\sigma=-s}^{s}(2s+1)\begin{Bmatrix} l_0 & l & l' \\ l_0 & s & l_0\end{Bmatrix}\begin{Bmatrix} l_0 & l'+l_0 & l' \\ l_0 & s & l_0\end{Bmatrix}\\
&\quad\times\sum_{m_{01},m_{02},m_{03},m_{04}}\begin{pmatrix} l_0 & l_0 & s \\ m_{01} & m_{02} & \sigma\end{pmatrix}\begin{pmatrix} l_0 & l_0 & s \\ m_{03} & m_{04} & \sigma\end{pmatrix}\\
&\qquad\qquad\times\frac{a_{l_0m_{01}}a_{l_0m_{02}}a_{l_0m_{03}}a_{l_0m_{04}}}{\widehat{C}^2_{l_0}},
\end{aligned}
$$

in view of (54) in the Appendix. Now the triangle conditions entail that the summands are nonzero only for $s\le 2l_0$; from (50) and (55) below we learn that

$$
\left|\begin{pmatrix} l_0 & l_0 & s \\ m_{01} & m_{02} & \sigma\end{pmatrix}\right|,\left|\begin{pmatrix} l_0 & l_0 & s \\ m_{03} & m_{04} & \sigma\end{pmatrix}\right|\le\frac{1}{\sqrt{2l_0+1}},
$$

$$
\left|\begin{Bmatrix} l_0 & l & l' \\ l_0 & s & l_0\end{Bmatrix}\right|,\left|\begin{Bmatrix} l_0 & l'+l_0 & l' \\ l_0 & s & l_0\end{Bmatrix}\right|\le\frac{1}{\sqrt{2l'+1}}\le\frac{1}{\sqrt{2l+1}},
$$

whereas it is also trivial to note that

$$
\left|\frac{a_{l_0m_{01}}a_{l_0m_{02}}}{\widehat{C}_{l_0}}\right|\le\max_m\frac{|a_{l_0m}|^2}{\widehat{C}_{l_0}}=(2l_0+1)\max_m\frac{|a_{l_0m}|^2}{\sum_{k=-l_0}^{l_0}|a_{l_0k}|^2}\le 2l_0+1.
$$

Hence

$$
\begin{aligned}
|B| &\leq \sum_{s\leq 2l_0}\sum_{|\sigma|\leq 2s}\frac{(2l_0+1)}{2l+1}\sum_{m_{01},m_{02},m_{03},m_{04}}\frac{1}{2l_0+1}\times(2l_0+1)^2 \\
&\leq \frac{(4l_0+1)(2l_0+1)^7}{2l+1}.
\end{aligned} \tag{39}
$$

Combining (36), (37), (38) and (39), we have

$$
|E\{X_{l,L}X_{l',L}|\Im_{[Lr]}\}| \leq \frac{1}{L}\frac{(4l_0+1)(2l_0+1)^7}{2l+1}\delta_{l+l_0}^{l'},
$$

$$
\left|\sum_{l=[Ls]}^{[Lt]}\sum_{l'=l+1}^{[Lt]}E\{X_{l,L}X_{l',L}|\Im_{[Lr]}\}\right| \leq \frac{\log L}{L}\frac{(4l_0+1)(2l_0+1)^7}{2l+1}.
$$

The remaining part of the proof is similar. More precisely, we have

$$
\begin{aligned}
&E\{X_{l,L}^2|\Im_{[Lr]}\} \\
&\quad = \frac{1}{2}E\{\widehat{I}_{l_0,l,l+l_0}^4-1|\Im_{[Lr]}\}, \\
&E\{\hat{I}_{l_0,l,l+l_0}^4|\Im_{[Lr]}\} \\
&\quad = 3\sum_{m_{01},m_{02},m_{03},m_{04}}\sum_{m_{11},m_{12},m_{21},m_{22}}\begin{pmatrix} l_0 & l & l+l_0 \\ m_{01} & m_{11} & m_{21}\end{pmatrix} \\
&\qquad\qquad \times\begin{pmatrix} l_0 & l & l+l_0 \\ m_{02} & -m_{11} & -m_{21}\end{pmatrix} \\
&\qquad \times\begin{pmatrix} l_0 & l & l+l_0 \\ m_{03} & m_{12} & m_{22}\end{pmatrix}\begin{pmatrix} l_0 & l & l+l_0 \\ m_{04} & -m_{12} & -m_{22}\end{pmatrix} \\
&\qquad \times\frac{a_{l_0m_{01}}a_{l_0m_{02}}a_{l_0m_{03}}a_{l_0m_{04}}}{\widehat{C}_{l_0}^2} \\
&\qquad + 6\sum_{m_{01},m_{02},m_{03},m_{04}}\sum_{m_{11},m_{12},m_{21},m_{22}}\begin{pmatrix} l_0 & l & l+l_0 \\ m_{01} & m_{11} & m_{21}\end{pmatrix} \\
&\qquad\qquad \times\begin{pmatrix} l_0 & l & l+l_0 \\ m_{02} & -m_{11} & m_{22}\end{pmatrix} \\
&\qquad \times\begin{pmatrix} l_0 & l & l+l_0 \\ m_{03} & m_{12} & -m_{21}\end{pmatrix}\begin{pmatrix} l_0 & l & l+l_0 \\ m_{04} & -m_{12} & -m_{22}\end{pmatrix} \\
&\qquad \times\frac{a_{l_0m_{01}}a_{l_0m_{02}}a_{l_0m_{03}}a_{l_0m_{04}}}{\widehat{C}_{l_0}^2} \\
&\quad = 3A+6B.
\end{aligned}
$$

As before,

$$A = \sum_{m_{01},m_{02},m_{03},m_{04}} \frac{\delta_{m_{01}}^{-m_{02}} \delta_{m_{03}}^{-m_{04}}}{(2l_0+1)^2} \frac{a_{l_0 m_{01}} a_{l_0 m_{02}} a_{l_0 m_{03}} a_{l_0 m_{04}}}{\widehat{C}_{l_0}^2}$$

$$= \sum_{m_{01},m_{03}} \frac{1}{(2l_0+1)^2} \frac{|a_{l_0 m_{01}}|^2 |a_{l_0 m_{03}}|^2}{\widehat{C}_{l_0}^2} = 1$$

and [see (54)]

$$\sum_{m_{11},m_{12},m_{21},m_{22}} \begin{pmatrix} l_0 & l & l+l_0 \\ m_{01} & m_{11} & m_{21} \end{pmatrix} \begin{pmatrix} l_0 & l & l+l_0 \\ m_{02} & -m_{11} & m_{22} \end{pmatrix}$$

$$\times \begin{pmatrix} l_0 & l & l+l_0 \\ m_{03} & m_{12} & -m_{21} \end{pmatrix} \begin{pmatrix} l_0 & l & l+l_0 \\ m_{04} & -m_{12} & -m_{22} \end{pmatrix}$$

$$= \sum_{s \le 2l_0} \sum_{\sigma=-s}^{s} (2s+1) \begin{Bmatrix} l & l+l_0 & l_0 \\ l_0 & s & l \end{Bmatrix}^2$$

$$\times \begin{pmatrix} l_0 & l_0 & s \\ m_{01} & m_{02} & \sigma \end{pmatrix} \begin{pmatrix} l_0 & l_0 & s \\ m_{03} & m_{04} & \sigma \end{pmatrix},$$

whence

$$B \le \frac{(4l_0+1)(2l_0+1)^7}{2l+1},$$

as in (39). We have thus shown that, for some $C > 0$,

$$\left| E\left\{ \left( \sum_{l=[Ls]}^{[Lt]} X_{l,L} \right)^2 \middle| \Im_{[Lr]} \right\} - (t-s) \right|$$

$$\le \left| \sum_{l=[Ls]}^{[Lt]} (E\{X_{l,L}^2 | \Im_{[Lr]}\} - 1) \right| + 2 \left| \sum_{l=[Ls]}^{[Lt]} \sum_{l'=l+1}^{[Lt]} E\{X_{l,L} X_{l',L} | \Im_{[Lr]}\} \right| + \frac{2}{L}$$

$$\le \frac{C}{L} \left\{ \left| \sum_{l=[Ls]}^{[Lt]} \frac{(4l_0+1)^8}{2l+1} \right| + \log L + 1 \right\}$$

$$\le C(4l_0+1)^8 \frac{\log L}{L} \to 0 \qquad \text{as } L \to \infty.$$

Thus the proof is complete. □

Theorem 5.1 can be immediately applied to derive the asymptotic distribution under the null of several non-Gaussianity tests. For instance, we might focus on $\sup_{0\le r\le 1} J_{aL;l_0,K}(r)$; by the continuous mapping theorem we obtain

$$\lim_{L\to\infty} P\Big\{\sup_{0\le r\le 1} J_{aL;l_0,K}(r)\le x\Big\} = P\Big\{\sup_{0\le r\le 1} W(r)\le x\Big\} = 2\Phi(x)-1, \qquad x\ge 0, a=1,2,3,4, \tag{40}$$

where $\Phi(\cdot)$ denotes the cumulative distribution function of a standard Gaussian variable (for the last equality, see, e.g., [5]). We shall now discuss the behavior of these procedures under some examples of non-Gaussian spherical fields. The behavior of higher-order angular power spectra under non-Gaussian alternatives is an extremely important research topic in modern cosmology, and still almost completely open for mathematical research. Very few analytic results are available, whereas the cosmological debate is still open on the nature of the non-Gaussianity to be expected. A simple and popular model for non-Gaussian temperature fluctuations reads

$$T_{\mathrm{NG}}(\theta,\varphi) = T(\theta,\varphi) + f_{\mathrm{NL}}\{T^2(\theta,\varphi) - ET^2(\theta,\varphi)\}; \tag{41}$$

as before, we take $T(\theta,\varphi)$ to be an isotropic Gaussian field with zero mean. For "small" values of the nonlinearity parameter $f_{\mathrm{NL}}$, (41) can be viewed as a general approximation for random fields with minor departures from Gaussianity: the quadratic term can be regarded as the leading factor in a Taylor expansion of a general field $g(T(\theta,\varphi))$, for a suitably regular function $g(\cdot)$. Equivalently, the terms on the left-hand side can be considered to be the first two elements of expansion of $g(T(\theta,\varphi))$ into a series of orthogonal Hermite polynomials $H_q(T)$, $q=1,2,\ldots.$ For these reasons, (41) is very widely adopted to represent the primordial field of temperature fluctuations in cosmological models of inflations, which stand now as the leading models for the dynamics in the primordial epochs around the Big Bang [3, 19, 28]. It is known in the physics literature [19] that the bispectrum of (41) can be approximated by

$$B_{l_1l_2l_3} = Gf_{\mathrm{NL}}h_{l_1l_2l_3}\begin{pmatrix} l_1 & l_2 & l_3\\ 0 & 0 & 0\end{pmatrix}\{C_{l_1}C_{l_2}+C_{l_2}C_{l_3}+C_{l_1}C_{l_3}\}, \tag{42}$$

where $G$ is a positive constant,

$$h_{l_1l_2l_3} = \left(\frac{(2l_1+1)(2l_2+1)(2l_3+1)}{4\pi}\right)^{1/2},$$

and lower-order terms are neglected. Expression (42) is known as the Sachs–Wolfe bispectrum; we shall take (42) as a benchmark model for our discussion of non-Gaussianity. We introduce a very mild regularity condition on the angular power spectrum, that is, we shall assume that $C_l$ is such that, for fixed $l_0>0$,

$$C_{l+l_0}\propto C_l, \tag{43}$$

where $\propto$ denotes that the ratio of the right- and left-hand sides tends to a positive constant, as $l \to \infty$. The assumption described by (43) is not unreasonable: indeed, $C_l \propto l^{-\alpha}$ (for some positive constant $\alpha$) for many, if not most, cosmological models.

For simplicity, let us assume that the normalizing angular power spectrum is nonrandom, that is, known a priori; without loss of generality we take $K = 0$. Recall that ([34], equations 8.1.2.12 and 8.5.2.32)

$$\begin{pmatrix} l_1 & l_2 & l_3 \\ 0 & 0 & 0 \end{pmatrix} = \frac{(-1)^{(l_1+l_2+l_3)/2}[(l_1+l_2+l_3)/2]!}{[(l_1+l_2-l_3)/2]![(l_1-l_2+l_3)/2]![(-l_1+l_2+l_3)/2]!} \times \left\{ \frac{(l_1+l_2-l_3)!(l_1-l_2+l_3)!(-l_1+l_2+l_3)!}{(l_1+l_2+l_3+1)!} \right\}^{1/2}.$$

Thus, for fixed $l_0 \geq 2$,

$$\begin{aligned} \begin{pmatrix} l_0 & l & l+l_0 \\ 0 & 0 & 0 \end{pmatrix} &= \frac{(-1)^{l_0+l}(l+l_0)\times\cdots\times(l+1)}{l_0!} \\ &\quad \times \frac{\sqrt{(2l_0)!}}{\sqrt{(2l+1)\times\cdots\times(2l+2l_0+1)}} \\ &= C\frac{(-1)^{l_0+l}}{\sqrt{l}} + O\left(\frac{1}{l^{3/2}}\right), \end{aligned}$$

for some $C > 0$ which depends on $l_0$ but not on $l$. Then we have easily that

$$\begin{aligned} EJ_{3L}(r) &\propto \frac{f_{\mathrm{NL}}}{\sqrt{L}} \sum_{l=l_0+1}^{[Lr]} \sqrt{\frac{(2l_0+1)(2l+1)(2l+2l_0+1)}{4\pi}} \\ &\quad \times \frac{1}{\sqrt{l}}\left[\sqrt{\frac{C_{l_0}C_l}{C_{l+l_0}}} + \sqrt{\frac{C_{l_0}C_{l+l_0}}{C_l}} + \sqrt{\frac{C_lC_{l+l_0}}{C_{l_0}}}\right] \\ &\propto \frac{f_{\mathrm{NL}}}{\sqrt{L}} \sum_{l=l_0+1}^{[Lr]} \sqrt{l}\sqrt{\frac{C_{l_0}C_l}{C_{l+l_0}}} \propto \frac{f_{\mathrm{NL}}}{\sqrt{L}} \sum_{l=l_0+1}^{[Lr]} \sqrt{l} \propto f_{\mathrm{NL}}L. \end{aligned} \tag{44}$$

Likewise

$$\begin{aligned} EJ_{4L}(r) &\propto \frac{1}{\sqrt{L}} \sum_{l=l_0+1}^{[Lr]} \{EI^2_{l_0,l,l+l_0} - 1\} \geq \frac{1}{\sqrt{L}} \sum_{l=1}^{[Lr]} \{(EI_{l_0,l,l+l_0})^2 - 1\} \\ &\propto \frac{f^2_{\mathrm{NL}}}{\sqrt{L}} \sum_{l=l_0+1}^{[Lr]} l \propto f^2_{\mathrm{NL}} L^{3/2}. \end{aligned} \tag{45}$$

Of course, both (44) and (45) diverge as the number of observed multipoles increases ($L \to \infty$), that is, as the resolution of the experiment improves. The constants of proportionality are typically small; for the model we adopt in the simulations below, they are of the order $10^{-4}$ for (44) and $10^{-8}$ for (45).

On the other hand, for $l_1 = l_2 = l_3 = l$,

$$\begin{pmatrix} l & l & l \\ 0 & 0 & 0 \end{pmatrix} = (-1)^{3l/2} \frac{[3l/2]!}{[(l/2)!]^3} \left\{ \frac{[l!]^3}{(3l+1)!} \right\}^{1/2}$$
$$\propto (-1)^{3l/2} \frac{[3l]^{(3l+1)/2}}{l^{3/2(l+1)}} \left\{ \frac{l^{3l+3/2}}{(3l)^{3(l+1/2)}} \right\}^{1/2} \propto \frac{(-1)^{3l/2}}{l},$$

where we have used Stirling's formula $n!/(\sqrt{2\pi} n^{n+1/2} e^{-n}) = 1 + O(12n^{-1})$. Now recall that

$$ET^2(\theta, \varphi) = \sum_{l=1}^{\infty} \frac{(2l+1)}{4\pi} C_l < \infty,$$

which implies $lC_l = o(l^{-1})$ (in realistic experiments, a weighting factor $G_l$ should be included to account for the beam pattern of the antenna). Hence we have

$$\begin{aligned} EJ_{1L;l_0,K}(r) &= O\left( \frac{1}{\sqrt{L}} \sum_{l=l_0+K}^{[Lr]} \sqrt{\frac{(2l+1)^3}{4\pi}} \begin{pmatrix} l & l & l \\ 0 & 0 & 0 \end{pmatrix} \frac{C_l^2}{\sqrt{C_l^3}} \right) \\ &= O\left( \frac{1}{\sqrt{L}} \sum_{l=l_0+K}^{[Lr]} \sqrt{lC_l} \right) = o(1) \qquad \text{as } L \to \infty. \end{aligned} \tag{46}$$

A similar heuristic argument can be used for $J_{2L;l_0,K}(r)$, suggesting that these statistics may have very little, if any, power against alternatives of the form (42).

To validate the previous heuristics, we present a small Monte Carlo study on the power of these testing procedures. To this aim, we generated 200 spherical Gaussian fields according to the currently favored scenario for CMB fluctuations, the so-called $\Lambda CDM$ model; we omit the details, which can be found, for instance, in [6]. Nonlinearities were then introduced by means of (41), which represents the simplest non-Gaussian model, as argued earlier. To ease comparisons with existing procedures, we follow the standard parametrization used in the astrophysical literature, yielding a variance for $T$ of the order of $\mathrm{Var}(T) = ET^2 \simeq 10^{-8}$. We can then present a rough relationship between the value of the nonlinearity parameter $f_{\mathrm{NL}}$ and the relative amount of the non-Gaussian signal, namely

$$\frac{\sqrt{\mathrm{Var}\{f_{\mathrm{NL}} T^2\}}}{\sqrt{\mathrm{Var}\{T\}}} = \sqrt{2} f_{\mathrm{NL}} \sqrt{\mathrm{Var}\{T\}} \simeq f_{\mathrm{NL}} \times 10^{-4}. \tag{47}$$

We consider $f_{\mathrm{NL}} = 0, 100, 300, 1000$ and we focus on the statistics $\sup_{0 \le r \le 1} J_{3L;l_0,K}(r)$, $\sup_{0 \le r \le 1} J_{4L;l_0,K}(r)$ for $l_0 = 2$ and $K = 0, 2, 4$; see (40).

We omit reporting results for $\sup_{0\le r\le 1} J_{1L;l_0,K}(r)$, $\sup_{0\le r\le 1} J_{2L;l_0,K}(r)$ because the power related to these procedures turned out to be negligible, as expected from (46). Other types of statistics may be considered without altering the main conclusions we are going to draw here. We start from the empirical sizes (type I errors), which are reported in Table 1; we write $S_{a,L}$ for $\sup_{0\le r\le 1} J_{a,L;2,K}(r)$, $a = 3, 4$. We take $L = 250, 500$; these values are conservative: we recall that $L$ is reckoned to be of the order of 600/800 for WMAP and 2000/2500 for Planck. Note that all values in Table 1 and those in Tables 3–5 are expressed as percentages (%).

Results in Table 1 suggest that the asymptotic theory presented in Theorem 5.1 provides a good approximation for the finite-sample behavior in the case of $J_{3,L}$; the approximation is slightly less satisfactory for $J_{4,L}$, but the results improve markedly with $L$, which is reassuring. Because the test statistics are free of nuisance parameters, it is also possible to derive directly threshold values under the null of Gaussianity by Monte Carlo replications. We adopt both approaches to derive the power properties reported in Tables 3–5; the Monte Carlo critical values are reported in Table 2.

Note how the tabulated values approach the asymptotic results [i.e., 1.645, 1.96 —see (40)] when $L$ increases. To analyze the power of the test, we consider $f_{\mathrm{NL}} = 100, 300, 1000$; from (47) we can argue heuristically that these values correspond approximately to 1%, 3% and 10% of non-Gaussianity in the maps, respectively. In Tables 3–5, T and A denote the power with respect to tabulated and asymptotic critical values, respectively.

At first sight, the results reported seem quite encouraging, as compared, for instance, to existing methods such as the empirical process, wavelets and local curvature: see, for instance, [6] for numerical simulations on the performance of these procedures. We stress once again, however, that this comparison could be to some extent misleading, insofar as in this paper we are sticking to some simplifying assumptions that are unrealistic for CMB experiments (namely, the absence of gaps and noise in the observed maps). The results we report, however, certainly suggest that further investigation of these procedures under more realistic circumstances is worth pursuing. We also note that the statistics based on $J_{3L;l_0,K}(r)$ substantially outperform those based on $J_{4L;l_0,K}(r)$ for this range of parameter values; indeed, for the latter the power is nonnegligible only when $f_{\mathrm{NL}} = 1000$, that is, when the level of non-Gaussianity in the map is approximately 10%. However, a comparison of (44) and (45) suggests that the power for $J_{4L;l_0,K}(r)$ may improve more rapidly as the resolution of the experiments grows; this is an important factor to keep in mind, as it is expected that the satellite Planck will achieve observations with $L$ of the order of 2500. Also, because $J_{4L;l_0,K}(r)$ does not depend on the sign of the sample bispectrum, it can be expected to be more robust against other types of non-Gaussian behavior. It is also worth noting that the power of the tests grows with the pooling parameter $K$; for values larger than 6–8, however, the effect is nearly negligible and the computational cost becomes prohibitive.

Table 1
*Empirical sizes* ($f_{\mathrm{NL}} = 0$), $\alpha = 10\%$ (5%)

| | $K = 0$ | $K = 2$ | $K = 4$ |
|---|---|---|---|
| $S_{3,250}$ | 9.5 (4.5) | 9.0 (4.0) | 8.0 (2.5) |
| $S_{3,500}$ | 10.5 (4.5) | 9.5 (4.5) | 9.5 (4.0) |
| $S_{4,250}$ | 5.5 (3.5) | 3.0 (1.5) | 1.0 (0.5) |
| $S_{4,500}$ | 6.5 (4.5) | 8.0 (4.0) | 6.0 (2.0) |

Table 2
*Monte Carlo critical values,* $\alpha = 10\%$ (5%)

| | $K = 0$ | $K = 2$ | $K = 4$ |
|---|---|---|---|
| $S_{3,250}$ | 1.61 (1.81) | 1.55 (1.83) | 1.61 (1.72) |
| $S_{3,500}$ | 1.69 (1.90) | 1.63 (1.92) | 1.61 (1.80) |
| $S_{4,250}$ | 1.44 (1.71) | 1.27 (1.51) | 1.06 (1.35) |
| $S_{4,500}$ | 1.56 (1.88) | 1.52 (1.85) | 1.38 (1.67) |

Table 3
*Rejection rates for* $f_{\mathrm{NL}} = 100$, $\alpha = 10\%$ (5%)

| | $K = 0$ (T) | $K = 2$ (T) | $K = 4$ (T) | $K = 0$ (A) | $K = 2$ (A) | $K = 4$ (A) |
|---|---|---|---|---|---|---|
| $S_{3,250}$ | 18.0 (10.5) | 20.0 (12.0) | 18.5 (17.5) | 17.5 (7.0) | 16.5 (9.5) | 17.5 (11.0) |
| $S_{3,500}$ | 20.5 (13.5) | 30.5 (20.0) | 39.5 (33.0) | 20.5 (12.5) | 30.5 (20.0) | 38.5 (26.0) |
| $S_{4,250}$ | 9.0 (5.0) | 10.5 (4.5) | 10.0 (4.5) | 5.5 (3.0) | 3.0 (1.5) | 1.0 (1.0) |
| $S_{4,500}$ | 10.0 (5.0) | 10.0 (6.0) | 10.0 (6.0) | 7.0 (4.0) | 7.0 (4.0) | 6.0 (2.0) |

Table 4
*Rejection rates for* $f_{\mathrm{NL}} = 300$, $\alpha = 10\%$ (5%)

| | $K = 0$ (T) | $K = 2$ (T) | $K = 4$ (T) | $K = 0$ (A) | $K = 2$ (A) | $K = 4$ (A) |
|---|---|---|---|---|---|---|
| $S_{3,250}$ | 31.0 (26.0) | 47.5 (38.5) | 56.5 (51.5) | 31.0 (21.5) | 43.0 (33.0) | 56.0 (41.5) |
| $S_{3,500}$ | 47.5 (43.0) | 88.5 (69.0) | 88.5 (86.5) | 48.5 (40) | 88.5 (68.5) | 88.5 (79.0) |
| $S_{4,250}$ | 10.5 (5.0) | 11.5 (5.0) | 9.0 (5.5) | 6.0 (3.0) | 3.5 (1.5) | 2.0 (0.5) |
| $S_{4,500}$ | 12.0 (6.0) | 11.5 (7.5) | 11.5 (7.0) | 9.0 (4.0) | 9.5 (6.0) | 7.5 (2.5) |

Table 5
*Rejection rates for* $f_{\mathrm{NL}} = 1000$, $\alpha = 10\%$ (5%)

| | $K = 0$ (T) | $K = 2$ (T) | $K = 4$ (T) | $K = 0$ (A) | $K = 2$ (A) | $K = 4$ (A) |
|---|---|---|---|---|---|---|
| $S_{3,250}$ | 80.0 (73.0) | 99.5 (99.5) | 100 (100) | 80.0 (70.5) | 99.5 (99.5) | 100 (99.5) |
| $S_{3,500}$ | 98.0 (98.0) | 100 (100) | 100 (100) | 98.0 (97.5) | 100 (100) | 100 (100) |
| $S_{4,250}$ | 15.5 (10.5) | 26.0 (18.5) | 32.5 (18.5) | 12.0 (7.5) | 10.5 (5.5) | 8.5 (3.5) |
| $S_{4,500}$ | 37.5 (29.0) | 57.5 (50.5) | 71.5 (58.5) | 34.0 (24.0) | 54.0 (47.0) | 60.5 (47.5) |

In our opinion, there are two remarkable features that emerge from this section. Equations (44) and (45) suggest that a testing procedure in harmonic space can yield consistent tests of Gaussianity even for fixed-radius, nonergodic fields, at least under the simplifying assumption of this paper. This is to some extent an unexpected result. The second remarkable fact is the huge impact of the choice of combined angular scales on the expected power under non-Gaussian alternatives. It is noteworthy that the common choice of a (close to) "main diagonal" configuration can yield negligible power, the expected value of the non-Gaussian signal decreasing to zero as the resolution of the experiment improves. The determination of the triples of angular scales $(l_1, l_2, l_3)$ where the largest part of the non-Gaussian signal is to be expected, for a given class of models, represents an issue of great importance for future cosmological data analysis.

**6. Comments and conclusion.** It is important to make clear that the asymptotic theory presented in this paper is of a rather different nature with respect to what is usually undertaken for random processes or fields. More precisely, we are not assuming that the information grows in the sense that a larger interval or region of observations becomes available, but rather we assume that the same region (spherical surface, in our case) is observed with greater and greater resolution. We labeled this framework high-resolution asymptotics, whereas the more standard case where the observed volume grows can be termed, as usual, large-sample asymptotics. The idea that some consistent inference can be drawn from a process observed on a fixed finite region is certainly not new; see, for instance, [32] for a fixed-domain asymptotics approach to study optimal linear prediction of spatial processes (kriging). The term infill asymptotics is also used in the statistical literature with the same meaning. We believe that this paradigm will become more and more fruitful in the years to come, with many possible contexts of applications. In the case of cosmological research, a proper understanding of the nature of the asymptotic theory involved is likely to be quite relevant from both the theoretical and the practical point of view. In fact, note that sequential experiments to measure CMB radiation result in exactly the same *last scattering* surface being measured, while the resolution improves steadily over time. For instance, the above mentioned NASA satellite WMAP, launched in 2001, is observing the same surface as the ESA mission Planck, due to be launched in 2007: in terms of the standard large-sample asymptotics, no improvement should be expected. On the other hand, for statistical properties that can be consistently investigated as the resolution of the experiment grows, Planck does offer substantial new information, its expected resolution outperforming WMAP by a factor 3 or 4 (in any case,

Planck will offer other improvements besides better angular resolution, e.g., better polarization measurements). It seems therefore quite relevant to suggest, as we did in this paper, that Gaussianity tests may exist that are high-resolution consistent, at least under some simplifying assumptions. The fact that consistent inferences can be drawn from nonergodic random fields defined on a bounded domain has other important consequences if we focus on epistemological issues. The status of cosmology as a science is occasionally questioned, on the grounds that it is, in a sense, a discipline based by definition on a single observation (our Universe). The possibility to draw consistent inferences for fixed-radius random fields provides, in our view, a strong argument to consider the corresponding physical properties fully within the domain of scientific investigation. It is a challenging task to characterize, under general conditions, the complete set of properties on which high-resolution consistent inferences can be drawn.

## APPENDIX

**A.1. The properties of Wigner's coefficients.** In this paper, we make extensive use of Wigner's $3j$ coefficients, which are a very powerful tool to represent properties of random fields which are invariant to rotations. These coefficients were introduced in the framework of the quantum theory of angular momenta; they are also widely used in algebra in the framework of representation theory. In this section we shall recall some of their properties for convenience.

Wigner's coefficients are defined implicitly by

$$\int_0^{\pi}\int_0^{2\pi} Y_{l_1m_1}(\theta,\varphi)Y_{l_2m_2}(\theta,\varphi)Y_{l_3m_3}(\theta,\varphi)\sin\theta\, d\varphi\, d\theta$$
$$= \left(\frac{(2l_1+1)(2l_2+1)(2l_3+1)}{4\pi}\right)^{1/2} \begin{pmatrix} l_1 & l_2 & l_3 \\ 0 & 0 & 0 \end{pmatrix} \begin{pmatrix} l_1 & l_2 & l_3 \\ m_1 & m_2 & m_3 \end{pmatrix}.$$

Many explicit representations are also available, but for general values of $l_i, m_i$ they are lengthy and hardly informative. For instance, it can be shown that ([34], expression 8.2.1.5)

$$\begin{pmatrix} l_1 & l_2 & l_3 \\ m_1 & m_2 & m_3 \end{pmatrix}$$
$$= (-1)^{l_3+m_3+l_2+m_2} \left[\frac{(l_1+l_2-l_3)!(l_1-l_2+l_3)!(l_1-l_2+l_3)!}{(l_1+l_2+l_3+1)!}\right]^{1/2}$$
$$\times \left[\frac{(l_3+m_3)!(l_3-m_3)!}{(l_1+m_1)!(l_1-m_1)!(l_2+m_2)!(l_2-m_2)!}\right]^{1/2}$$
$$\times \sum_z \frac{(-1)^z(l_2+l_3+m_1-z)!(l_1-m_1+z)!}{z!(l_2+l_3-l_1-z)!(l_3+m_3-z)!(l_1-l_2-m_3+z)!},$$

where the summation runs over all $z$'s such that the factorials are nonnegative. We list here some important properties, and refer to [34] for proofs and further discussion:

(a) Wigner's $3j$ coefficients are real valued;
(b) they are different from zero only if $m_1+m_2+m_3=0$;
(c) (parity) for any triple $l_1,l_2,l_3$

$$\begin{pmatrix} l_1 & l_2 & l_3 \\ m_1 & m_2 & m_3 \end{pmatrix} = (-1)^{l_1+l_2+l_3} \begin{pmatrix} l_1 & l_2 & l_3 \\ -m_1 & -m_2 & -m_3 \end{pmatrix};$$

(d) (symmetry) for any triple $l_1,l_2,l_3$

$$\begin{aligned} \begin{pmatrix} l_1 & l_2 & l_3 \\ m_1 & m_2 & m_3 \end{pmatrix} &= \begin{pmatrix} l_2 & l_3 & l_1 \\ m_2 & m_3 & m_1 \end{pmatrix} = \begin{pmatrix} l_3 & l_1 & l_2 \\ m_3 & m_1 & m_2 \end{pmatrix} \\ &= (-1)^{l_1+l_2+l_3} \begin{pmatrix} l_3 & l_2 & l_1 \\ m_3 & m_2 & m_1 \end{pmatrix} \\ &= (-1)^{l_1+l_2+l_3} \begin{pmatrix} l_1 & l_3 & l_2 \\ m_1 & m_3 & m_2 \end{pmatrix} \\ &= (-1)^{l_1+l_2+l_3} \begin{pmatrix} l_2 & l_1 & l_3 \\ m_2 & m_1 & m_3 \end{pmatrix}; \end{aligned}$$

(e) (orthonormality) for any triple $l_1,l_2,l_3$

$$\sum_{m_1=-l_1}^{l_1} \sum_{m_2=-l_2}^{l_2} \sum_{m_3=-l_3}^{l_3} \begin{pmatrix} l_1 & l_2 & l_3 \\ m_1 & m_2 & m_3 \end{pmatrix}^2 = 1 \tag{48}$$

and

$$\sum_{m_1=-l_1}^{l_1} \sum_{m_2=-l_2}^{l_2} \begin{pmatrix} l_1 & l_2 & L \\ m_1 & m_2 & M \end{pmatrix} \begin{pmatrix} l_1 & l_2 & L' \\ m_1 & m_2 & M' \end{pmatrix} = \frac{\delta_L^{L'} \delta_M^{M'}}{2L+1}; \tag{49}$$

(f) (upper bound) for any $l_1,l_2,l_3$

$$\begin{pmatrix} l_1 & l_2 & l_3 \\ m_1 & m_2 & m_3 \end{pmatrix} = O([\max\{l_1,l_2,l_3\}]^{-1/2}); \tag{50}$$

(g1) (sum of coefficients, I) for any positive integers $a,b$

$$\sum_\alpha (-1)^{-\alpha} \begin{pmatrix} a & a & b \\ \alpha & -\alpha & \beta \end{pmatrix} = (-1)^a \sqrt{2a+1}\delta_b^0 \delta_\beta^0; \tag{51}$$

(g2) (sums of coefficients, II) for any positive integers $a,b,c,d,e$ and $f$

$$\begin{aligned} &\sum_{\alpha,\beta,\gamma} \sum_{\varepsilon,\delta,\phi} (-1)^{e+f+\varepsilon+\phi} \begin{pmatrix} a & b & e \\ \alpha & \beta & \varepsilon \end{pmatrix} \begin{pmatrix} c & d & e \\ \gamma & \delta & -\varepsilon \end{pmatrix} \\ &\qquad \times \begin{pmatrix} a & d & f \\ \alpha & \delta & -\phi \end{pmatrix} \begin{pmatrix} c & b & f \\ \gamma & \beta & \phi \end{pmatrix} \\ &= \begin{Bmatrix} a & b & e \\ c & d & f \end{Bmatrix}; \end{aligned} \tag{52}$$

(g3) (sums of coefficients, III) for any positive integers $a, b, c, d, e$ and $f$ (see [34], 8.7.3.12)

$$\sum_{\alpha,\beta,\delta} (-1)^{\delta+\gamma} \begin{pmatrix} a & b & c \\ -\alpha & -\beta & \gamma \end{pmatrix} \begin{pmatrix} a & f & d \\ \alpha & \varphi & -\delta \end{pmatrix} \begin{pmatrix} e & b & d \\ -\varepsilon & \beta & -\delta \end{pmatrix} = (-1)^{b-f} \begin{pmatrix} c & f & e \\ \gamma & \varphi & -\varepsilon \end{pmatrix} \begin{Bmatrix} a & b & c \\ e & f & d \end{Bmatrix}; \tag{53}$$

(g4) (sums of coefficients, IV) for any positive integers $a, b, c, d, e$ and $f$ (see [34], 8.7.4.20)

$$\begin{aligned} &\sum_{\beta,\gamma,\varepsilon,\varphi} \begin{pmatrix} a & b & c \\ \alpha & -\beta & -\gamma \end{pmatrix} \begin{pmatrix} d & f & e \\ \delta & -\varphi & -\varepsilon \end{pmatrix} \begin{pmatrix} g & b & e \\ \eta & \beta & \varepsilon \end{pmatrix} \begin{pmatrix} j & f & c \\ \mu & \varphi & \gamma \end{pmatrix} \\ &\quad = (-1)^{a-b+c+d+e-f} \\ &\qquad \times \sum_{s\sigma} \begin{pmatrix} a & s & j \\ \alpha & \sigma & -\mu \end{pmatrix} \begin{pmatrix} g & s & d \\ \eta & \sigma & -\delta \end{pmatrix} \begin{Bmatrix} b & c & a \\ j & s & f \end{Bmatrix} \begin{Bmatrix} b & e & g \\ d & s & f \end{Bmatrix}. \end{aligned} \tag{54}$$

Equation (52) can be used as the definition of *Wigner's 6j coefficient*, which appears on the right-hand side (with curly brackets). We refer again to ([34], Chapter 9) for some (extremely complicated) explicit expression for these coefficients and many of their properties. For our purposes it is sufficient to recall the following:

(h) for any positive integers $a, b, c, d, e$ and $f$

$$\left| \begin{Bmatrix} a & b & c \\ d & e & f \end{Bmatrix} \right| \leq \min\left( \frac{1}{\sqrt{(2c+1)(2f+1)}}, \frac{1}{\sqrt{(2a+1)(2d+1)}}, \frac{1}{\sqrt{(2b+1)(2e+1)}} \right); \tag{55}$$

(i) the $6j$ coefficient is invariant under any permutation of its columns.

### A.2. Proofs of technical lemmas.

Proof of Proposition 3.1. Let $\{I_1, \ldots, I_p\}$ be a partition of $I$ into $2 \times 3$ matrices made up with two of its rows; write $\mathcal{I}$ for the class of these partitions, which has cardinality $(2p-1)!!$ (the number of combinations by which we can match $2p$ rows two by two to form $p$ pairs). It suffices to notice that

$$D[\Gamma_P(I,3); l_1, l_2, l_3] = \sum_{\{I_1,\ldots,I_p\} \in \mathcal{I}} \prod_{i=1}^{p} D[\Gamma_P(I_i,3); l_1, l_2, l_3].$$

For any $I_i$, $D[\Gamma_P(I_i,3);l_1,l_2,l_3]$ is easily seen to include $\Delta_{l_1l_2l_3}$ nonzero summands, each of them of the same form, up to a relabeling of the indexes; more precisely,

$$\begin{aligned} &D[\Gamma_P(I_i,3);l_1,l_2,l_3] \\ &\quad = \Delta_{l_1l_2l_3} \sum_{m_{i1},m_{i3},m_{i3}} (-1)^{m_{i1}+m_{i2}+m_{i3}} \begin{pmatrix} l_1 & l_2 & l_3 \\ m_{i1} & m_{i2} & m_{i3} \end{pmatrix} \\ &\qquad \times \begin{pmatrix} l_1 & l_2 & l_3 \\ -m_{i1} & -m_{i2} & -m_{i3} \end{pmatrix} \\ &\quad = \Delta_{l_1l_2l_3} \sum_{m_{i1},m_{i3},m_{i3}} \begin{pmatrix} l_1 & l_2 & l_3 \\ m_{i1} & m_{i2} & m_{i3} \end{pmatrix}^2 = \Delta_{l_1l_2l_3}. \end{aligned}$$

Thus

$$D[\Gamma_P(I,3);l_1,l_2,l_3] = \sum_{\{I_1,\dots,I_p\}\in I} \prod_{i=1}^{p} \Delta_{l_1l_2l_3} = (2p-1)!!\Delta_{l_1l_2l_3}^p,$$

as claimed. □

Proof of Proposition 3.2. (a) The case $\#(I)=4$. This result can be established as a straightforward application of Lemmas 3.1–3.3.

(b) The case $\#(I)=6$. Without loss of generality we can take $I=\{1,2,\dots,6\}$, and show that

$$\sum_{m_{11}=-l_1}^{l_1} \cdots \sum_{m_{63}=-l_3}^{l_3} \prod_{i=1}^{6} \begin{pmatrix} l_1 & l_2 & l_3 \\ m_{i1} & m_{i2} & m_{i3} \end{pmatrix} \sum_{\gamma\in\Gamma(6,3)\backslash\Gamma_P(6,3)} \delta(\gamma;l_1,l_2,l_3) = O(l_1^{-1}).$$

It is sufficient to consider only the diagrams with no flat edges; it is readily seen that the latter can be partitioned into (a) the unconnected (unpaired) diagrams and (b) the connected diagrams. Now for (a) we note that, because there cannot be any flat edge, if the diagram is unconnected but not paired the set of rows must necessarily be partitioned into a group of two and a group of four; in other words, after some rearrangement of indices it must be possible to write any diagram $\gamma$ belonging to (a) as $\gamma=\gamma_1\cup\gamma_2$, where $\gamma_1\in\Gamma_P(2,3)$ and $\gamma_2\in\Gamma_C(4,3)$. Thus, exactly as shown above, we obtain that the corresponding terms are bounded by $C(2l_1+1)^{-1}$. It suffices then to look at the connected diagrams. Consider first $\gamma\in\Gamma_{CL(2)}(6,3)$; from Lemma 3.2 and simple manipulations we obtain immediately

$$D[\Gamma_{CL(2)}(6,3)] = O(l_1^{-1}D[\Gamma_C(4,3)]) = O(l_1^{-2}). \tag{56}$$

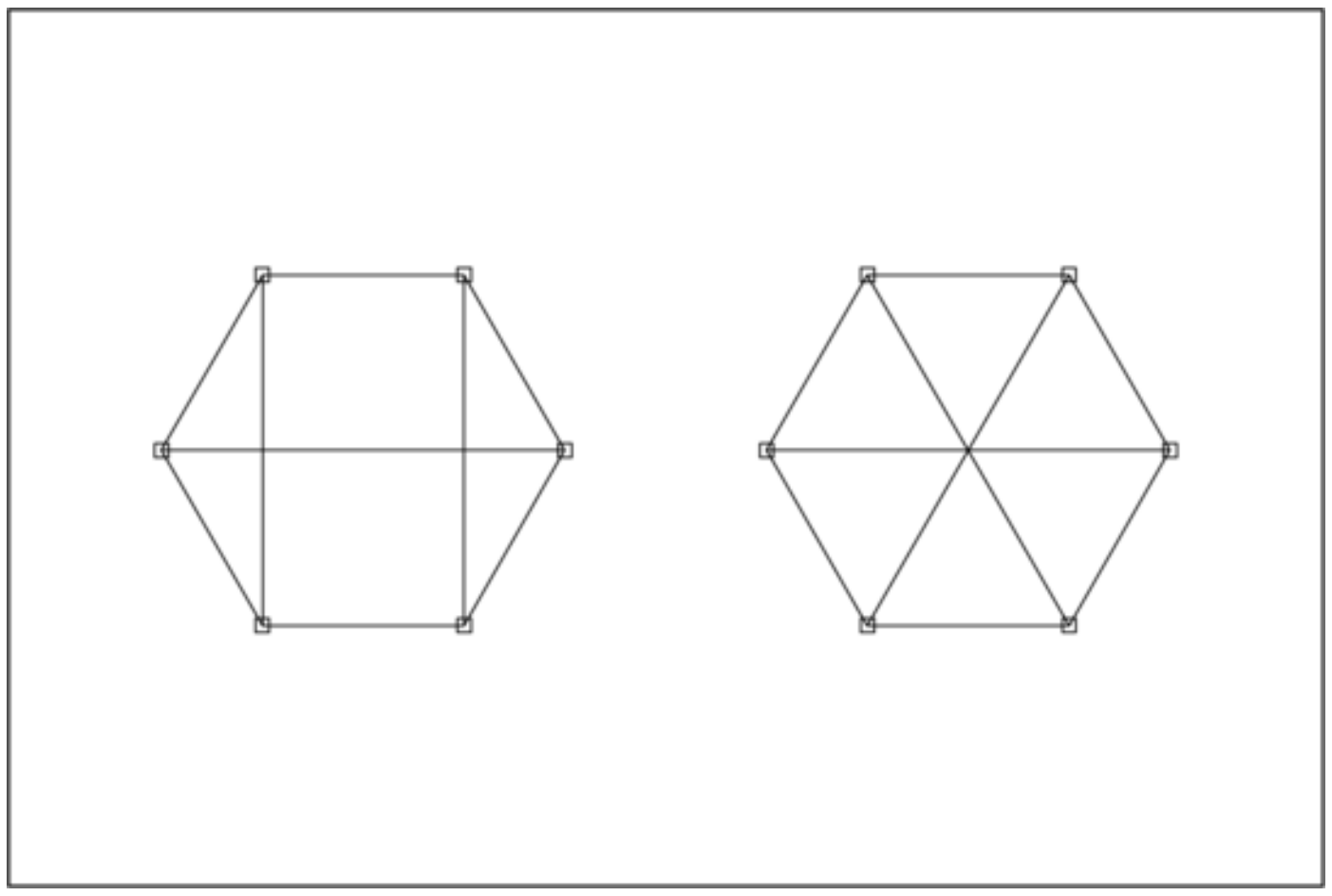

FIG. 4. $\gamma \in \Gamma_{CL(3)}(6,3)$, $\gamma \in \Gamma_{CL(4)}(6,3)$.

In case $\gamma \in \Gamma_{C\overline{L(2)}}(6,3)$, it can be readily shown that we must have $\gamma \in \{\Gamma_{CL(3)}(6,3) \cup \Gamma_{CL(4)}(6,3)\}$; in other words, these diagrams have at least a loop of order either 3 or 4 [$\Gamma_{C\overline{L(4)}}(6,3)$ is empty]. The latter claim can be established as follows. We have a graph with six vertices, each of which is of degree 3. Hence we can argue as in ([11], page 48) to show that the closure of the graph is complete (it is a clique). Then, from Theorem 4.5 of ([11], page 48) it follows that the graph is Hamiltonian, that is, it contains a spanning cycle. With a permutation of the indices, we can then take the vertices to be at the corners of a regular hexagon, and we are left with an edge free for each of them. These edges will connect different vertices to form loops of order 3 or 4. The graphs corresponding to $\gamma \in \Gamma_{CL(3)}(6,3)$ with at least a loop of order 3 are labeled type A; those with all loops of order at least 4 [$\gamma \in \Gamma_{CL(4)}(6,3)$] are labeled type B (see Figure 4).

Let $\gamma$ correspond to a type A graph; by Lemma 3.3 we easily have (see Figure 3)

$$D[\gamma; l_1, l_2, l_3] = O\left(\left\{\begin{pmatrix} l_1 & l_2 & l_3 \\ l_1 & l_2 & l_3 \end{pmatrix}\right\} D[\Gamma_{C\overline{L(2)}}(4,3); l_1, l_2, l_3]\right)$$

$$= O\left(\begin{Bmatrix} l_1 & l_2 & l_3 \\ l_1 & l_2 & l_3 \end{Bmatrix}^2\right) = O(l_3^{-2}).$$

For type B, again up to a permutation of the indices we can take, with no loss of generality,

$$m_{11} = -m_{21}, \qquad m_{12} = -m_{42}, \qquad m_{13} = -m_{63},$$
$$m_{22} = -m_{32}, \qquad m_{23} = -m_{53}, \qquad m_{31} = -m_{61},$$
$$m_{33} = -m_{43}, \qquad m_{41} = -m_{51}, \qquad m_{52} = -m_{62}.$$

The corresponding expected value is premultiplied by the factor

$$(-1)^{m_{11}+m_{12}+m_{13}+m_{31}+m_{32}+m_{33}+m_{51}+m_{52}+m_{53}} = 1,$$

whence we have

$$\begin{aligned}
&\left\{\prod_{i=1,3,5}\prod_{j=1}^{3}\sum_{m_{ij}}\right\}\begin{pmatrix} l_1 & l_2 & l_3 \\ m_{11} & m_{12} & m_{13}\end{pmatrix} \\
&\qquad\qquad \times \begin{pmatrix} l_1 & l_2 & l_3 \\ -m_{11} & -m_{32} & -m_{53}\end{pmatrix}\begin{pmatrix} l_1 & l_2 & l_3 \\ m_{31} & m_{32} & m_{33}\end{pmatrix} \\
&\qquad\qquad \times \begin{pmatrix} l_1 & l_2 & l_3 \\ -m_{51} & -m_{12} & -m_{33}\end{pmatrix}\begin{pmatrix} l_1 & l_2 & l_3 \\ m_{51} & m_{52} & m_{53}\end{pmatrix} \\
&\qquad\qquad \times \begin{pmatrix} l_1 & l_2 & l_3 \\ -m_{31} & -m_{52} & -m_{13}\end{pmatrix} \\
&\quad = \sum_x (-1)^{2x}(2x+1)\begin{Bmatrix} l_1 & l_2 & l_3 \\ l_3 & l_3 & x\end{Bmatrix}\begin{Bmatrix} l_1 & l_2 & l_3 \\ l_2 & x & l_2\end{Bmatrix}\begin{Bmatrix} l_1 & l_2 & l_3 \\ x & l_1 & l_1\end{Bmatrix},
\end{aligned} \tag{57}$$

where we have used ([34], equation 10.2.3.17, page 339, and equation 10.2.4.20, page 340). In view of (55) we easily obtain

$$(57) \le \frac{C}{\sqrt{l_1 l_2 l_3}}\sum_{x=1}^{l_3}\frac{1}{x} = O(l_1^{-1/2} l_2^{-1/2} l_3^{-1/2}\log l_3) = O(l_1^{-1}).$$

(c) The case $\#(I) = 8$. Again, we can take $I = \{1, 2, \ldots, 8\}$ and show that

$$\begin{aligned}
&\sum_{m_{11}=-l_1}^{l_1}\cdots\sum_{m_{83}=-l_3}^{l_3}\prod_{i=1}^{8}\begin{pmatrix} l_1 & l_2 & l_3 \\ m_{i1} & m_{i2} & m_{i3}\end{pmatrix} \\
&\qquad \times \sum_{\gamma\in\{\Gamma(8,3)\setminus\Gamma_P(8,3)\}} \delta(\gamma; l_1, l_2, l_3) = O(l_1^{-1}).
\end{aligned}$$

As before, it is readily seen that the diagrams with no flat edges can be partitioned into (a) the unconnected unpaired diagrams and (b) the connected diagrams. Now for (a) we note that

$$\begin{aligned}
&\gamma \in [\{\Gamma_{\overline{C}}(8,3)\setminus\Gamma_P(8,3)\} \cap \Gamma_{\overline{F}}(8,3)] \\
&\quad \Rightarrow \gamma \in [\{\Gamma_C(I_1,3)\otimes\Gamma_C(I_2,3)\} \cup \{\Gamma_C(I_3,3)\otimes\Gamma_C(I_4,3)\} \\
&\qquad\qquad \cup \{\Gamma_C(I_5,3)\otimes\Gamma_C(I_6,3)\otimes\Gamma_C(I_7,3)\}],
\end{aligned}$$

where $(I_1, I_2), (I_3, I_4), (I_5, I_6, I_7)$ are partitions of $\{1, 2, \ldots, 8\}$ into disjoint sets such that

$$\begin{aligned}
&\#(I_1) = 6, \qquad \#(I_2) = 2, \\
&\#(I_3) = 4, \qquad \#(I_4) = 4, \\
&\#(I_5) = 4, \qquad \#(I_6) = 2, \qquad \#(I_7) = 2;
\end{aligned}$$

by $\Gamma_1 \otimes \Gamma_2$, we denote the set of all diagrams of the form $\gamma = \gamma_1 \cup \gamma_2$, for $\gamma_1 \in \Gamma_1$ and $\gamma_2 \in \Gamma_2$. For any two families of diagrams $\Gamma_1(I_1, 3), \Gamma_2(I_2, 3)$ such that $I_1 \cap I_2 = \varnothing$, it is readily checked that

$$D[\Gamma_1 \otimes \Gamma_2; l_1, l_2, l_3] = O(D[\Gamma_1; l_1, l_2, l_3] \times D[\Gamma_2; l_1, l_2, l_3]),$$
$$D[\Gamma_1 \cup \Gamma_2; l_1, l_2, l_3] = O(D[\Gamma_1; l_1, l_2, l_3] + D[\Gamma_2; l_1, l_2, l_3]);$$

thus

$$\begin{aligned} &D[\Gamma_{\overline{C}}(8,3) \backslash \Gamma_P(8,3); l_1, l_2, l_3] \\ &\quad = O(D[\Gamma_C(6,3) \otimes \Gamma_C(2,3); l_1, l_2, l_3]) \\ &\qquad + O(D[\Gamma_C(4,3) \otimes \Gamma_C(4,3); l_1, l_2, l_3]) \\ &\qquad + O(D[\Gamma_C(4,3) \otimes \Gamma_C(2,3) \otimes \Gamma_C(2,3); l_1, l_2, l_3]) \\ &\quad = O(D[\Gamma_C(6,3); l_1, l_2, l_3]) \\ &\qquad + O(D^2[\Gamma_C(4,3); l_1, l_2, l_3]) + O(D[\Gamma_C(4,3); l_1, l_2, l_3]) \\ &\quad = O(l_1^{-1}), \end{aligned}$$

as shown previously. It suffices then to look at the connected diagrams. Diagrams with a 2-loop can be handled by Lemma 3.2, and then by using results on lower-order diagrams. So we just have to consider $\gamma \in \Gamma_{C\overline{L(2)}}(8,3)$; these diagrams must have at least a loop of order 3 or 4 [in other words, $\Gamma_{C\overline{L(4)}}(8,3)$ is empty]. The latter claim can be established in the same manner as before; we repeat the argument for completeness. We have a graph with eight vertices, each of which is of degree 3. Hence we can argue as in ([11], page 48) to show that the closure of the graph is complete. Then, from Theorem 4.5 of [11] it follows that the graph is Hamiltonian, that is, it contains a spanning cycle. With a permutation of the indices, we can then take the vertices to be at the corners of a regular octagon, and we are left with an edge free for each of them. These edges will connect different vertices to form loops of order 3, 4 or 5. Graphs with loops of order 3 can be dealt with as before through Lemma 3.3. If all loops are of even order, then the graph is bipartite ([11], Theorem 2.4, page 23); we label it a type C graph, for which we provide two isomorphic representations in Figure 5.

If there is at least a loop of order 5, we label it a type D graph, for which two isomorphic representations are provided in Figure 6.

To analyze the behavior of the components $D[\gamma, l_1, l_2, l_3]$ corresponding to type C graphs, we can take with no loss of generality

$$\begin{aligned} &m_{11} = -m_{21}, &&m_{12} = -m_{42}, &&m_{13} = -m_{83}, &&m_{32} = -m_{22}, \\ &m_{31} = -m_{41}, &&m_{33} = -m_{63}, &&m_{51} = -m_{61}, &&m_{52} = -m_{82}, \\ &m_{53} = -m_{43}, &&m_{71} = -m_{81}, &&m_{72} = -m_{62}, &&m_{73} = -m_{23}, \end{aligned}$$

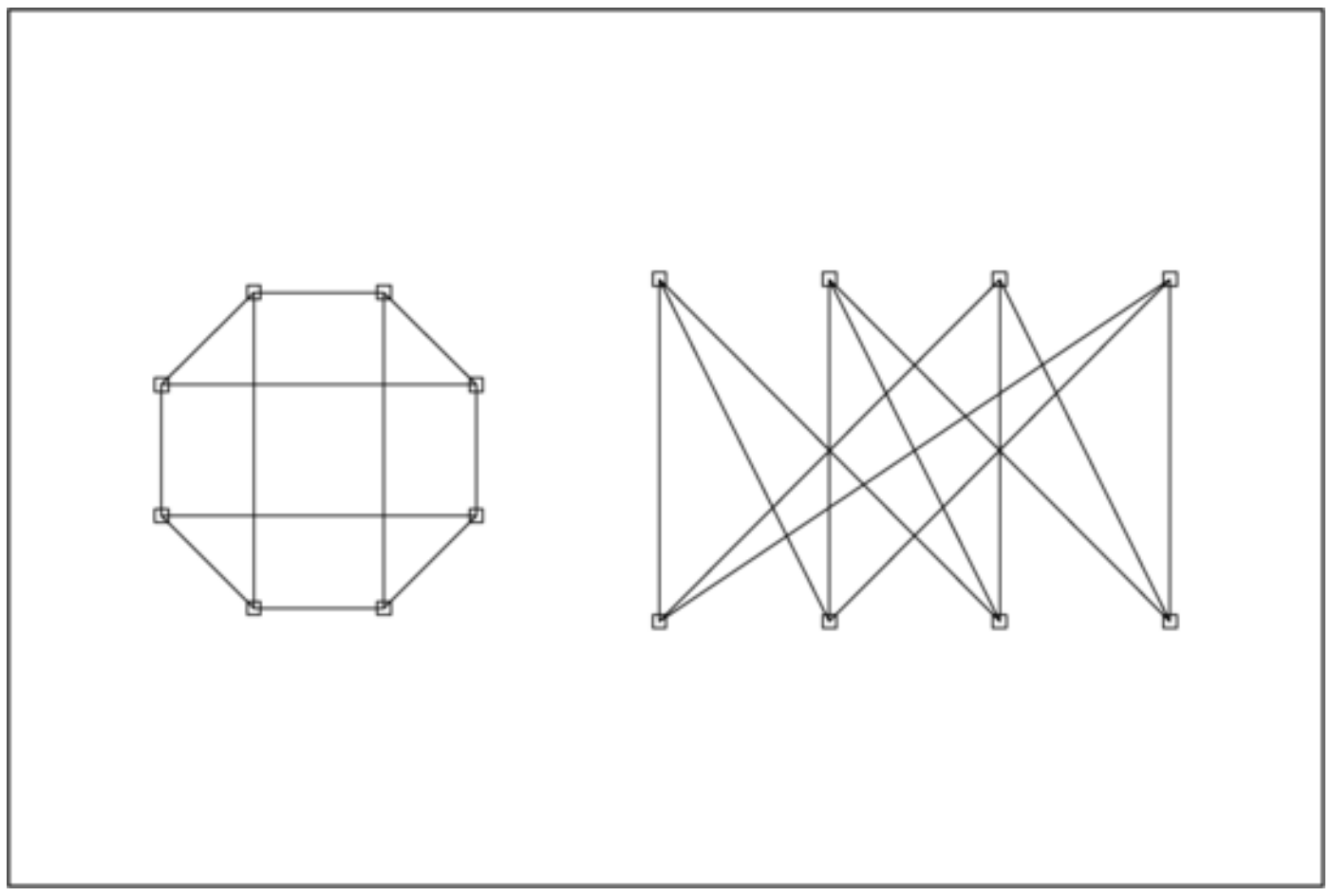

FIG. 5. *Isomorphic type* C *graphs.*

which leads to

$$
\begin{aligned}
&\left\{\prod_{i=1,3,5,7}\prod_{j=1}^{3}\sum_{m_{ij}}\right\}\begin{pmatrix} l_1 & l_2 & l_3 \\ m_{11} & m_{12} & m_{13}\end{pmatrix}\begin{pmatrix} l_1 & l_2 & l_3 \\ -m_{11} & -m_{32} & -m_{73}\end{pmatrix} \\
&\qquad\times\begin{pmatrix} l_1 & l_2 & l_3 \\ m_{31} & m_{32} & m_{33}\end{pmatrix}\begin{pmatrix} l_1 & l_2 & l_3 \\ -m_{31} & -m_{12} & -m_{53}\end{pmatrix} \\
&\qquad\times\begin{pmatrix} l_1 & l_2 & l_3 \\ m_{51} & m_{52} & m_{53}\end{pmatrix} \\
&\qquad\times\begin{pmatrix} l_1 & l_2 & l_3 \\ -m_{51} & -m_{72} & -m_{33}\end{pmatrix}\begin{pmatrix} l_1 & l_2 & l_3 \\ m_{71} & m_{72} & m_{73}\end{pmatrix} \\
&\qquad\times\begin{pmatrix} l_1 & l_2 & l_3 \\ -m_{71} & -m_{52} & -m_{13}\end{pmatrix} \\
&= (-1)^{l_1-l_3-l_1+l_3}\sum_x (2x+1)\begin{Bmatrix} l_2 & l_1 & x \\ l_3 & l_3 & l_2\end{Bmatrix}\begin{Bmatrix} l_2 & l_1 & x \\ l_3 & l_1 & l_3\end{Bmatrix} \\
&\qquad\times\begin{Bmatrix} l_3 & l_3 & x \\ l_2 & l_2 & l_1\end{Bmatrix}\begin{Bmatrix} l_3 & l_1 & x \\ l_2 & l_2 & l_3\end{Bmatrix},
\end{aligned} \tag{58}
$$

where we have used ([34], equations 10.13.3.23 and 10.13.3.25, page 367); the sum runs over all positive integers $x$ which satisfy the triangle inequalities $l_2 - l_1 \le x \le l_2 + l_1$. By using (55), we obtain

$$
\begin{aligned}
(58) &\le C(2l_3+1)^2 \max_x \begin{Bmatrix} l_2 & l_1 & x \\ l_3 & l_3 & l_2\end{Bmatrix}\begin{Bmatrix} l_2 & l_1 & x \\ l_3 & l_1 & l_3\end{Bmatrix}\begin{Bmatrix} l_3 & l_3 & x \\ l_2 & l_2 & l_1\end{Bmatrix}\begin{Bmatrix} l_3 & l_1 & x \\ l_2 & l_2 & l_3\end{Bmatrix} \\
&= O(l_2^{-3/2} l_1^{-1/2}).
\end{aligned}
$$

In view of ([34], equations 10.13.1.1 and 10.13.1.3, page 361) the proof can be completed by an analogous argument for components corresponding to type D graphs. □

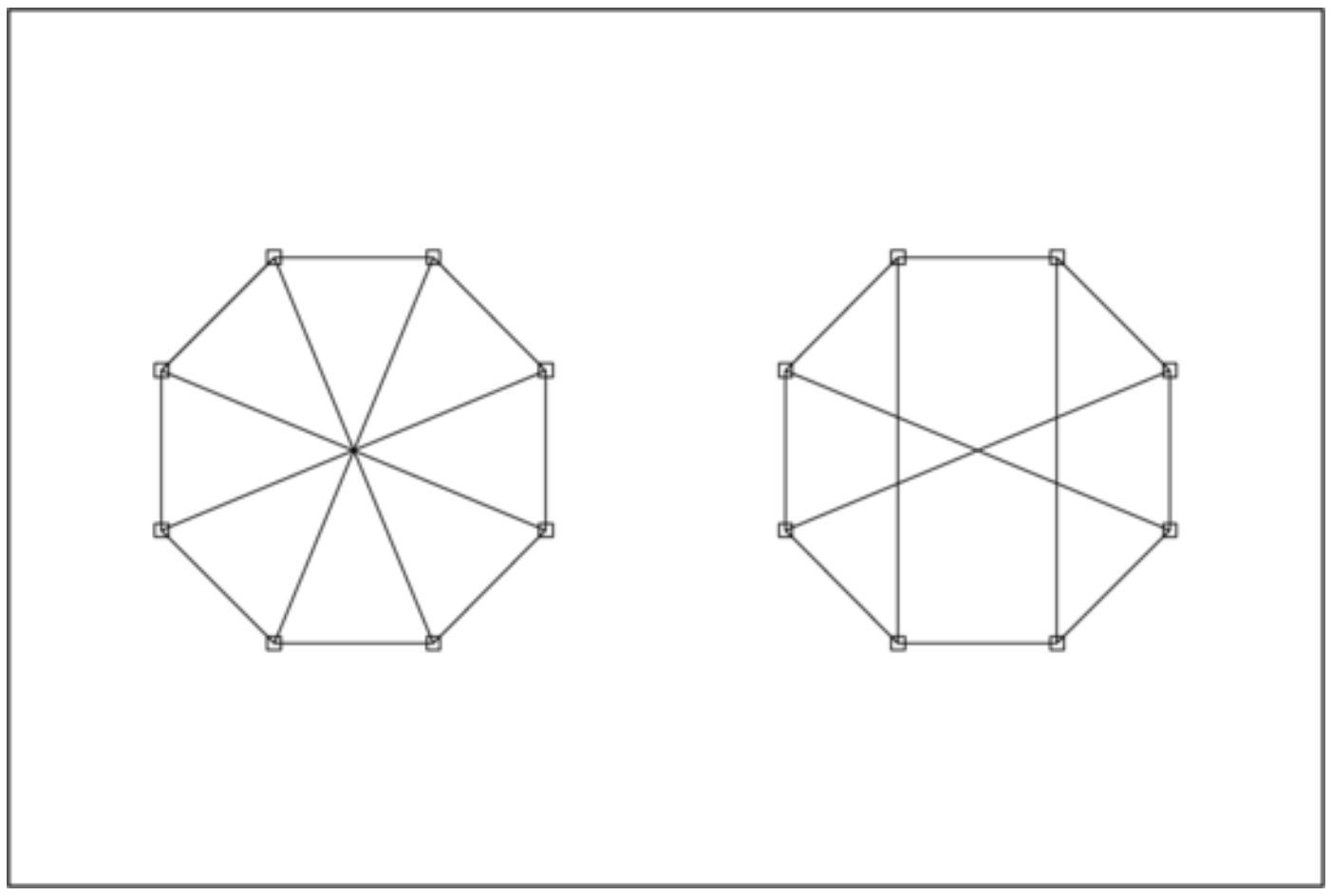

Fig. 6. *Isomorphic type* D *graphs.*

**Acknowledgments.** I am grateful to an Associate Editor and three anonymous referees for valuable comments that greatly improved the presentation of the paper; I am also grateful to P. Cabella for carrying out the simulations in Section 5.

Dipartimento di Matematica
Università di Roma "Tor Vergata"
via della Ricerca Scientifica 1
00133 Roma
Italy
E-mail: marinucc@mat.uniroma2.it